\documentstyle[12pt]{article}
\makeatletter
\newtheorem{lemma}{Lemma}[section]
\newcommand{\Lemma}{\begin{lemma}}
\newcommand{\EndLemma}{\end{lemma}}
\newtheorem{prop}[lemma]{Proposition}
\newcommand{\Prop}{\begin{prop}}
\newcommand{\EndProp}{\end{prop}}
\newtheorem{remark}[lemma]{\it Remark\/}
\newcommand{\Remark}{\begin{remark}}
\newcommand{\EndRemark}{\end{remark}}
\newtheorem{definition}[lemma]{\it Definition\/}
\newcommand{\Definition}{\begin{definition}}
\newcommand{\EndDefinition}{\end{definition}}
\newtheorem{example}[lemma]{\it Example\/}
\newcommand{\Example}{\begin{example}}
\newcommand{\EndExample}{\end{example}}
\newtheorem{note}[lemma]{\it Note\/}
\newcommand{\Note}{\begin{note}}
\newcommand{\EndNote}{\end{note}}
\newtheorem{thm}[lemma]{Theorem}
\newcommand{\Theorem}{\begin{thm}}
\newcommand{\EndTheorem}{\end{thm}}
\newtheorem{cor}[lemma]{Corollary}
\newcommand{\Cor}{\begin{cor}}
\newcommand{\EndCor}{\end{cor}}
\newcommand{\Proof}{\smallbreak\noindent{\it Proof.}\ \ }
\newcommand{\EndProof}{\QED}
\@addtoreset{equation}{section}
\makeatother
\catcode`\@=11
\font\tenmsy=msbm10 scaled 1200
\font\sevenmsy=msbm9
\font\fivemsy=msbm7

\newfam\msyfam
\textfont\msyfam=\tenmsy \scriptfont\msyfam=\sevenmsy
\scriptscriptfont\msyfam=\fivemsy
\def\Bbb{\ifmmode\let\next\Bbb@\else
 \def\next{\errmessage{Use \string\Bbb\space only in math mode}}\fi\next}
\def\Bbb@#1{{\Bbb@@{#1}}}
\def\Bbb@@#1{\fam\msyfam#1}
\def\hexnumber@#1{\ifnum#1<10 \number#1\else
 \ifnum#1=10 A\else\ifnum#1=11 B\else\ifnum#1=12 C\else
 \ifnum#1=13 D\else\ifnum#1=14 E\else\ifnum#1=15 F\fi\fi\fi\fi\fi\fi\fi}
\def\msx@{\hexnumber@\msxfam}
\def\msy@{\hexnumber@\msyfam}
\catcode`\@=\active
\newcommand{\Equation}{\begin{equation}}
\newcommand{\EndEquation}{\end{equation}}
\def\EQ#1{{\def\normalbaselines{\baselineskip15pt \lineskip6pt \lineskiplimit4pt}\,\vcenter{\normalbaselines\mathsurround=0pt\ialign{\hfil$\displaystyle##$\hfil&&$\,$$\displaystyle##$\hfil\crcr\mathstrut\crcr\noalign{\kern-\baselineskip}#1\crcr\mathstrut\crcr\noalign{\kern-\baselineskip}}}\,}}

\def\QED{\ifhmode\unskip\nobreak\fi\hfill{\it QED.}\bigbreak\penalty0\noindent}
\def\barb#1{{\bf\bar{\fam=-1 #1}}}
\def\tildeb#1{{\bf\tilde{\fam=-1 #1}}}

\newcommand{\GL}{\hbox{\sl GL}}
\newcommand{\Aut}{{\rm Aut\,}}
\newcommand{\ZZ}{{\Bbb Z}}
\newcommand{\CC}{{\Bbb C}}
\newcommand{\RR}{{\Bbb R}}
\newcommand{\NN}{{\Bbb N}}
\newcommand{\Ker}{{\rm Ker}\,}
\newcommand{\End}{{\rm End}\,}
\newcommand{\Mvir}{{M}}
\newcommand{\Lvir}{{L}}
\def\basis#1{{[#1]}} 
\newcommand{\character}{{\rm ch}\,}
\newcommand{\casimir}{{\kappa}}
\newcommand{\Trace}{{\rm Tr}\,}
\newcommand{\idem}{{e}}
\newcommand{\R}{{R}}
\newcommand{\1}{{\bf 1}}

\newcommand{\compo}{{\raise1pt\hbox{$\,\scriptstyle \circ\,$}}}
\newcommand{\Cyc}{{\rm Cyc}}
\newcommand{\Sym}{{\rm Sym}}
\newcommand{\<}{{\langle}}
\renewcommand{\>}{{\rangle}}
\newcommand{\onehalf}{\raise2pt\hbox{$\scriptstyle\frac12$}}
\newcommand{\CS}{{\cal S}}
\begin{document}
\title{
{\small\null
\vskip-10ex
\begin{flushright}
%Version 1.0,\ \ July, 2000\\
\end{flushright}
}
\vskip5ex
{Norton's Trace Formulae for the Griess Algebra\\
of a Vertex Operator Algebra\\
with Larger Symmetry
\vskip2ex
}
\author{\cr
Atsushi MATSUO
\cr\cr
{\normalsize\sl Department of Pure Mathematics and Mathematical Statistics,}\cr
{\normalsize\sl Centre for Mathematical Sciences,}\cr
{\normalsize\sl University of Cambridge}\cr
{\normalsize\sl Wilberforce Road, Cambridge CB3 0WB, England$\;^\dagger$}
%\cr\cr
%{\normalsize and}
%\cr\cr
%{\normalsize\sl Graduate School of Mathematical Sciences}\cr
%{\normalsize\sl University of Tokyo, Komaba, Tokyo 153-8914, Japan}
}
\date{}
}
\maketitle
{{\catcode`\@=11\long\def\@makefntext#1{\setbox0=\hbox{$\dagger.\ $}\leftskip1em\noindent\kern-\wd0\noindent#1}\catcode`\@=\active\footnotetext{$\dagger.\ $On leave of absence from the {\sl Graduate School of Mathematical Sciences, University of Tokyo, Komaba 3-8-1, Tokyo 153-8914, Japan}, supported by the Overseas Research Scholarship of the Ministry of Education, Science, Sports and Culture, Japan.
}}}
\noindent
{\small {\bf Abstract. }
Formulae expressing the trace of the composition of several (up to
five) adjoint actions of elements of the Griess algebra of a vertex
operator algebra are derived under certain assumptions on the action
of the automorphism group.  They coincide, when applied to the
moonshine module $V^\natural$, with the trace formulae obtained in a
different way by S.~Norton, and the spectrum of some idempotents
related to 2A, 2B, 3A and 4A element of the Monster is determined by
the representation theory of Virasoro algebra at $c=1/2$, $W_3$
algebra at $c=4/5$ or $W_4$ algebra at $c=1$.  The generalization to
the trace function on the whole space is also given for the
composition of two adjoint actions, which can be used to compute the
McKay-Thompson series for a 2A involution of the Monster.  }
\section*{Introduction}
Since Griess' construction of the Monster simple group \cite{Griess1}
as the automorphism group of a commutative nonassociative algebra of
dimension $196883+1$, many attempts are made in order to understand
the nature of this algebra.  Conway \cite{Conway} reconstructed a
slightly modified version of the algebra, called the Conway-Griess
algebra, and gave a description of a 2A involution (a transposition)
in terms of the eigenspace decomposition with respect to the adjoint
action of an idempotent of a particular type called the transposition
axis.  Some more formulae related to the decomposition of this kind
are obtained by Norton \cite{Norton}.  In particular, he wrote down a
trace formula for the composition of several (up to five) adjoint
actions of elements of the algebra:
$$\Trace {\it ad}_{a},\
\Trace {\it ad}_a{\it ad}_b,\ldots.
$$
These results are based on the explicit construction of the algebra as
well as the character table of the Monster and some of its subgroups.

On the other hand, Frenkel et al.\ \cite{FLM1} constructed a graded
vector space $V^\natural$ called the moonshine module, and showed that
the Conway-Griess algebra is naturally realized as the subspace of
degree (conformal weight) $2$ of $V^\natural$.  The multiplication and
the inner product of the algebra are actually a part of an infinite
series of bilinear operations on the whole space $V^\natural$ giving
this space with the structure of a vertex operator algebra (VOA)
\cite{Borcherds1},\cite{FLM2}. The moonshine module $V^\natural$ with
the VOA structure is probably the most natural object to be considered
in the study of the Monster in its relation to the moonshine
phenomena, as it was more or less clear in the original construction
of $V^\natural$ in \cite{FLM1}, and further supported by Borcherds's
solution \cite{Borcherds2} to the Conway-Norton conjecture
\cite{ConwayNorton}.  The algebra formed by the degree $2$ subspace of
a VOA is generally called the Griess algebra of the VOA.

Recently an attempt to understand the Conway-Griess algebra from the
VOA point of view was made by Miyamoto \cite{Miyamoto1}, who opened a
way to study the action of Monster elements on the moonshine module by
using a subVOA whose fusion rules have a nice symmetry.  In
particular, 2A involution of the Monster is described as an
automorphism obtained by the eigenspace decomposition with respect to
the action of Virasoro $L_0$ corresponding to the subVOA isomorphic to
$\Lvir(\frac12,0)$, which reproduces, on the level of the Griess
algebra, the action described by Conway mentioned above.  Further, the
structure of $V^\natural$ as a module over the tensor product of $48$
copies of $\Lvir(\frac12,0)$, called a frame, was studied in
\cite{DongGriessHohn}, and the VOA structure of the moonshine module
was reconstructed from a frame in \cite{Miyamoto3}.  In particular,
the character and the McKay-Thompson series for 2A, 2B and 3C elements
of the Monster can be computed by using a frame.

\vskip1ex

The primary purpose of this paper is to derive Norton's trace formula
concerning the adjoint action of elements of the Griess algebra
mentioned above from the VOA structure on the moonshine module
$V^\natural$ without using any explicit structure of the Griess
algebra or the Monster; the only particular property necessary for our
derivation is the fact that the components of $V^\natural$ fixed by
the full automorphism group coincides with the Virasoro submodule
generated by the vacuum vector up to degree 11.  This property says
that, while the automorphism group is finite as it is the Monster, the
symmetry of the VOA is large enough to separate the trivial components
into the small subspace of which the action on the VOA is determined
by the action of the Virasoro algebra (up to that degree).

It should be emphasized that, in our derivation of the formulae, we do
not need to know that the automorphism group is the Monster so that
our method is generally applicable to any VOA with the same property;
We call it a VOA of class~$\CS^n$ if its trivial components with
respect to the full automorphism group coincides with the Virasoro
submodule generated by the vacuum vector up to degree $n$.  This is
what we mean by a VOA with larger symmetry.

We will actually show in this paper (Theorem \ref{Theorem1}) that
the trace of the composition of $m$ adjoint actions of elements of the
Griess algebra, $m=1,\ldots,5$, is expressed in the same way as
Norton's formula with coefficients being replaced by certain rational
functions of the rank (central charge) $c$ and the dimension $d$ of
the Griess algebra if the VOA is of class~$\CS^{2m}$ under some
technical assumptions.

In spite that the trace must be invariant under the cyclic permutation
of the order of the operators involved in the trace, the explicit
expression for $m=4$ and $5$ does not satisfy this property in
general; it means that there is a restriction on the pair $(c,d)$
coming from our assumptions.  Also, in a slightly different way, we
see that if the Griess algebra contains an idempotent with central
charge being different from $0$ and $c$ then similar restrictions are
imposed on the pair $(c,d)$.  In this way, we may list up the possible
pairs of $(c,d)$ of a VOA satisfying our assumptions (Section 3).  In
particular, if the VOA is of class $\CS^{8}$ and it has an idempotent
as above, then the rank must be $24$ and the dimension must be
$196884$, i.e., those of the moonshine module $V^\natural$ (Theorem
\ref{Theorem2}).

Now, as we have established Norton's trace formulae in a different
way, we may use them to study the action of some elements of the
Monster.  Indeed, the trace formulae has sufficient information to
determine the spectrum of some idempotents related to 2A, 2B, 3A and
4A element of the Monster by the representation theory of
$\Lvir(\frac12,0)$, $W_3$ algebra at $c=4/5$ or $W_4$ algebra at $c=1$
(Subsection 4.2), reproducing some of the results of Conway
\cite{Conway} and Norton \cite{Norton} in the opposite way.

We note that the trace formulae would be generalized to the traces on
the higher degree subspaces.  In fact, we will show (Theorem
\ref{Theorem3}) that, under suitable assumptions, the trace functions
$$
\Trace a_{(1)}q^{L_0},\quad \hbox{and}\quad \Trace a_{(1)}b_{(1)}q^{L_0},
$$
where $a,b$ are elements of the Griess algbera, are indeed expressed
in terms of the character $\character V$ and its derivatives with
coefficients written by $(a|\omega)$, $(b|\omega)$ and $(a|b)$ as well
as the Eisenstein series $E_2(q)$ and $E_4(q)$ by using some
identities established by Zhu \cite{Zhu}.  As a corollary, we show
that the McKay-Thompson series $T_{\rm 2A}(q)$ for a 2A involution of
the Monster can be computed from $\character V^\natural=q(J(q)-744)$
and the characters of $\Lvir(\frac12,h)$, $h=0,1/2,1/16$, without
using a frame.

We finally note that our consideration based on the nonexistence of a
Monster invariant primary vector of degree less than $12$ in the
moonshine module can be understood to be complementary to the result
of Dong and Mason \cite{DongMason} that the one-point function for a
Monster invariant primary vector of degree $12$ gives rise to the cusp
form $\Delta(q)$ of weight~$12$.

\vskip1ex

Most of the results of this paper were obtained by using computer. The
author used Mathematica Ver.\ 3.0 for Linux.
\section{Preliminaries}
In this section, we recall or give some definitions and facts
necessary in ths paper.
\subsection{The Griess algebra of a vertex operator algebra}
Let $V$ be a vertex operator algebra (VOA) over the field $\CC$ of
complex numbers.  Recall that it is a vector space equipped with a
linear map
$$Y(a,z)=\sum_{n\in\ZZ}a_{(n)}z^{-n-1}\in (\End V)[[z,z^{-1}]]$$
and nonzero vectors $\1$ and $\omega$ satisfying a number of conditions \cite{Borcherds1},\cite{FLM2}.
We recall some of the properties (cf.\ \cite{MatsuoNagatomo}).
The operators $a_{(n)}$, ($a\in V, n\in\ZZ$), are subject to
\Equation
\sum_{i=0}^\infty{p\choose i}(a_{(r+i)}b)_{(p+q-i)}
=\sum_{i=0}^\infty (-1)^r{r\choose i}\Bigl(a_{(p+r-i)}b_{(q+i)}-(-1)^rb_{(q+r-i)}a_{(p+i)}\Bigr).
\label{eq:borcherds}
\EndEquation
where $a,b,c\in V$ and $p,q,r\in\ZZ$.
The {\it vacuum vector}\/ $\1$ satisfies
\Equation
\1_{(n)}a=\left\{\matrix{0,&(n\ne -1),\hfill \cr a,&(n=-1),\cr}\right.%\}
\quad
\hbox{and}
\quad
a_{(n)}\1=\left\{\matrix{0,&(n\geq 0),\hfill \cr a,&(n=-1).\cr}\right.%\}
\label{EqId}\EndEquation
The {\it conformal vector}\/ $\omega$ generates a representation of the Virasoro algebra:
$$
[L_m,L_n]=(m-n)L_{m+n}+\frac{m^3-m}{12}\delta_{m+n,0}\,c
$$
where $L_m=\omega_{(m+1)}$. The central charge $c$ is called the {\it
rank}\/ of the VOA.  The operator $L_{-1}$ satisfies 
\Equation
(L_{-1}a)_{(n)}=(\omega_{(0)}a)_{(n)}=-na_{(n-1)}.
\label{EqDer}\EndEquation

The operator $L_0$ is supposed to be semisimple giving rise to a
grading $ V=\bigoplus_{n=0}^\infty V^n $ of the VOA $V$, where $V^n$
denotes the eigenspace with eigenvalue $n$, which is supposed to be
finite-dimensional by definition. The eigenvalue is called the {\it
degree}.  It follows that $V^i_{(n)}V^j\subset V^{i+j-n-1}$ and $\1\in
V^0$.  We denote the sum of the subspaces with degree up to $n$ as
\Equation
V^{\leq n}=\bigoplus_{m=0}^n V^m.
\EndEquation

Throughout the paper, we assume that the grading of the VOA $V$ is of the form
\Equation
V=\bigoplus_{n=0}^\infty V^n,\quad\hbox{where}\quad V^0=\CC\1 \quad\hbox{and}\quad V^1=0.
\label{EqDoublyBosonic}
\EndEquation
Consider the first nontrivial subspace $B=V^2$. Set
\Equation
ab=a_{(1)}b,\quad
(a|b)\1=a_{(3)}b,\quad\hbox{for $a,b\in B$.}
\EndEquation
Then the multiplication gives a commutative nonassociative algebra
structure on $B$ and $(\ |\ )$ is a symmetric invariant bilinear form
on it.  The space $B$ equipped with these structures is called the
{\it Griess algebra\/} of $V$.  We denote the adjoint action of an
element $a\in B$ as
\Equation
\R_a\,:\,B\rightarrow B,\ x\mapsto xa\ (=ax).  
\EndEquation

By a slight abuse of terminology, we call a vector $\idem\in B$
satisfying $\idem_{(1)}\idem=2\idem$ an {\it idempotent}.  A vector
$\idem\in B$ generates a representation of the Virasoro algebra on $V$
if and only if it is an idempotent of $B$, for which the central
charge is given by $2(\idem|\idem)$.  The conformal vector $\omega$ of
the VOA $V$ is twice an identity element of the algebra, i.e., $\omega
a=2a$ for any $a\in B$.  The squared norm $(\omega|\omega)=c/2$ is
half the rank of the VOA $V$.

Recall that the VOA $V$ carries a unique invariant bilinear form $(\
|\ )$ up to normalization (\cite{Li}).  We normalize the form by
$(\1|\1)=1$.  It is indeed an extension of the form $(\ |\ )$ on $B$
to the whole space $V$, so we have denoted it by the same symbol.  We
note that $(V^i|V^j)=0$ if $i\ne j$. 
The relation
\Equation (a_{(n)}u|v)=(u|a_{(2m-2-n)}v),\quad (u,v\in V),
\EndEquation 
for a vector $a\in V^m$ such that $L_1a=0$ is a particular case of the
invariance of the form.  Therefore, $(a_{(n)}u|v)=(u|a_{(2-n)}v)$
holds for any $a\in B$ thanks to the assumption
(\ref{EqDoublyBosonic}).
\subsection{Virasoro submodule generated by the vacuum vector}
Let $V$ be a VOA and let $V_\omega$ be the Virasoro submodule generated by
the vacuum vector $1$ with respect to the action $L_n$ associated with
the conformal vector $\omega$.  Then, since $L_n\1=0$ for $n\geq -1$,
we have the following sequence of surjective homomorphisms of Virasoro
modules:
\Equation
\Mvir(c,0)/\Mvir(c,1)\rightarrow V_\omega\rightarrow \Lvir(c,0),
\label{eq:maps}
\EndEquation
where $\Mvir(c,h)$ (resp.\ $\Lvir(c,h)$) denote the Verma module
(resp.\ irreducible module) over the Virasoro algebra of central
charge $c$ with highest weight (lowest conformal weight) $h$.  Let us
denote the highest weight vector of $\Mvir(c,0)/\Mvir(c,1)$ mapped to
$\1\in V_\omega\subset V$ by the same symbol $\1$.

Now, let $P_n$ denote the set of all partitions of $n$ by integers
greater than $1$:
\Equation
P_n=\{\vec{m}=m_1m_2\cdots m_k\,|\,k, m_1,\ldots,m_k\in\NN,\ m_1\geq
m_2\geq\cdots\geq m_k\geq 2\}.
\label{eq:partition}
\EndEquation
For each partition $\vec{m}\in P_n$, we set
\Equation
\basis{\vec{m}}=\basis{m_1,m_2,\ldots,m_k}=L_{-m_1}L_{-m_2}\cdots L_{-m_k}\1\in
\Mvir(c,0)/\Mvir(c,1).
\EndEquation
Then the set $\{\basis{\vec{m}}\,|\, \vec{m}\in P_n\}$ forms a basis
of the subspace with conformal weight $n$ of the module
$\Mvir(c,0)/\Mvir(c,1)$.  We denote
\Equation
L_{\vec{m}}=L_{m_k}\cdots L_{m_2}L_{m_1}
\EndEquation
for a partition $\vec{m}=m_1m_2\cdots m_{k}$.

Recall that a singular vector (or a primary vector) of a Virasoro
module is a nonzero vector $v$ of the module such that
\Equation
L_m v=0\quad\hbox{for all $m\geq 1$}.
\label{eq:sing}
\EndEquation
In this paper, by convention, a {\it singular vector}\/ means a
nonzero vector in a highest weight module over the Virasoro algebra
satisfying (\ref{eq:sing}) which is not a multiple of the highest
weight vector (vacuum vector), and a {\it primary vector}\/ means any
nonzero vector in a VOA satisfying (\ref{eq:sing}).

The module $\Mvir(c,0)/\Mvir(c,1)$ contains a singular vector of conformal
weight up to $n$ if and only if the central charge $c$ is a zero of a
certain reduced polynomial $D_n(c)$, which can be computed by the
Kac-determinant formula.  We normalize the polynomials for $n=2,\ldots,10$ as follows.
\Equation
\EQ{
&D_2(c)=c,\cr
&D_4(c)=c(5c+22),\cr
&D_6(c)=c(2c-1)(5c+22)(7c+68),\cr
&D_8(c)=c(2c-1)(3c+46)(5c+3)(5c+22)(7c+68),\cr
&D_{10}(c)=10c(2c-1)(3c+46)(5c+3)(5c+22)(7c+68)(11c+232).\cr
}
\label{eq:Dn}
\EndEquation

If $D_n(c)\ne0$, then, up to the degree $n$, the maps (\ref{eq:maps})
are isomorphisms and we may identify $V_\omega$ with
$\Mvir(c,0)/\Mvir(c,1)$,
\subsection{VOA of class $\CS^n$}
An automorphism of a VOA $V$ is a linear isomorphism $g:V\rightarrow
V$ satisfying $g(a_{(n)}b)=(ga)_{(n)}(gb)$ for all $a,b\in V$ and all
$n\in\ZZ$ that fixes the conformal vector $\omega$.  Let $\Aut V$
denote the group of all automorphisms of $V$.

Any automorphism sends the vacuum vector $\1$ to itself and preserves
the grading.  Since an automorphism $g\in \Aut V$ fixes any vector in
the subspace $V_\omega$, the group $\Aut V$ acts on the quotient space
$V/V_\omega$ and on its graded pieces $V^n/V_\omega^n$.

\Definition
\rm
A VOA $V$ is said to be {\it of class~$\CS^n$}\/ if the action of
$\Aut V$ on $V^{\leq n}/V_\omega^{\leq n}$ is fixed-point free.
\EndDefinition

In other words, a VOA $V$ is of class~$\CS^n$ if $V^{\leq n}$ has no
extra fixed-vector other than those belonging to $V_\omega^{\leq n}$.

The VOA $\Lvir(c,0)$ associated with the irreducible highest
weight representation of the Virasoro algebra at central charge $c$
(cf.\ \cite{FrenkelZhu}) is obviously of class~$\CS^n$.  However,
since its Griess algebra $B$ is one-dimensional, this example is not
of our interest, although we will use $\Lvir(\frac12,0)$ etc.\ later
in a different context.

The main example of our concern is the moonshine module $V^\natural$
constructed by Frenkel et al.\ in \cite{FLM1} \cite{FLM2}.  It follows
from \cite{ConwayNorton} and \cite{Borcherds2} (cf.\ \cite{HaradaLang}
and \cite{DongMason}) that the VOA $V^\natural$ is of
class~$\CS^{11}$.  We will later see that the moonshine module
$V^\natural$, for which $c=24$ and $\dim B=196884$, has exceptionally
large symmetry in our sense.

\vskip1ex

In the rest of this paper, we always assume that $D_n(c)\ne0$ whenever
$V$ is supposed to be of class~$\CS^n$, for the cases when $D_n(c)=0$
are not interesting from our point of view.  The VOA's of our concern
are at most of class $\CS^{10}$ so that the excluded values of $c$ are
only $0$, $1/2$, $-46/3$, $-3/5$, $-22/5$, $-68/7$ and $-232/11$.
\section{Trace Formulae for the Griess Algebra of a VOA}
Let $B$ be the Griess algebra of a VOA $V$ of rank $c$, and let $d$
denote the dimension of $B$.  

In the first subsection, we give formulae that describes the traces
$\Trace R_{a_1}R_{a_2}\cdots R_{a_m}$ up to $m=5$ under appropriate
assumptions.  Sketch of the derivation of the formulae will be given
in subsequent subsections.
\subsection{The formulae}
To describe the formulae, we set
\Equation
(a_1|a_2|a_3)=(a_1|a_2a_3),
\EndEquation
which is a totally symmetric trilinear form on $B$, and define a
totally antisymmetric quinery form on $B$ by setting
\Equation
(a_1,a_2,a_3,a_4,a_5)\1=\frac{1}{5!}\sum_{\sigma\in S_5}
(-1)^{\ell(\sigma)}\sigma(a_1{}_{(3)}a_2{}_{(2)}a_3{}_{(1)}a_4{}_{(0)}a_5).
\label{eq:quinery}
\EndEquation
Here we let $\sigma\in S_5$ act by the permutation of the indices of
$a_i$, $i=1,\ldots,5$.

Let $\Cyc$ denote the operation of summing over the cyclic
permutations of the indices, and $\Sym$ denote that over all
permutations such that the results are apparently distinct after
performing the symmetries $a_ia_j=a_ja_i$, $(a_i|a_j)=(a_j|a_i)$ and
$(a_ia_j|a_k)=(a_i|a_ja_k)$ for any $i,j,k=1,\ldots,5$; For instance,
$$
\EQ{
&\Sym\, (a_1|a_2)(a_3|\omega)(a_4|\omega)\cr
&\quad=
(a_1|a_2)(a_3|\omega)(a_4|\omega)
+(a_1|a_3)(a_2|\omega)(a_4|\omega)
+(a_1|a_4)(a_2|\omega)(a_3|\omega)\cr
&\qquad +(a_2|a_3)(a_1|\omega)(a_4|\omega)
+(a_2|a_4)(a_1|\omega)(a_3|\omega)
+(a_3|a_4)(a_1|\omega)(a_2|\omega).\cr
}
$$

The result is summarized in the following theorem.

\Theorem
\sl
Let $B$ be the Griess algebra of a VOA $V$ such that the bilinear form
$(\ |\ )$ on $B$ is nondegenerate or $\Aut V$ is finite\footnote{We
suppose this for simplicity although the condition can be slightly
weakened.}
\item
{\rm (1)}
If $V$ is of class~$\CS^2$ then, for any $a\in B$, 
$$\Trace \R_{a} =\frac{4d}{c}(a|\omega).$$
\item
{\rm (2)}
If $V$ is of class~$\CS^4$ then, for any $a_1,a_2\in B$, 
$$
\Trace \R_{a_1}\R_{a_2}=\frac{-2(5c^2-88d+2cd)}{c(5c+22)}(a_1|a_2)+\frac{4(5c+22d)}{c(5c+22)}(a_1|\omega)(a_2|\omega).
$$
\item
{\rm (3)}
If $V$ is of class~$\CS^6$ then, for any $a_1,a_2,a_3\in B$, 
$$\EQ{
\Trace& \R_{a_1}\R_{a_2}\R_{a_3}\cr
&=\frac{-3{c^2}(70{c^2}+769c-340)+2d(4{c^3}-445{c^2}+12236c-5984)}{c(2c-1)(5c+22)(7c+68)}(a_1|a_2|a_3)\cr
&\hfill+\frac{4c(70{c^2}+1017c-340)-8d(32{c^2}-1419c+748)}{c(2c-1)(5c+22)(7c+68)}\,\Cyc\,(a_1|a_2)(a_3|\omega)\cr
&\hfill+\frac{5952 c(d-1)}{c(2c-1)(5c+22)(7c+68)}(a_1|\omega)(a_2|\omega)(a_3|\omega).\cr
}
$$
\item
{\rm (4)} If $V$ is of class~$\CS^8$ then, for any $a_1,a_2,a_3,a_4\in B$,
$$\EQ{
\Trace& \R_{a_1}\R_{a_2}\R_{a_3}\R_{a_4}\cr
&=\frac{1}{D_8(c)}
\Bigl(A_{1}(a_1a_2|a_3a_4)
+A_{2}(a_1a_3|a_2a_4)
+A_{3}(a_1a_4|a_3a_2)\cr
&\hfill+B\,\Sym\,(a_1|a_2|a_3)(a_4|\omega)
+C\,\Sym\,(a_1|a_2)(a_3|a_4)\hfill\cr
&\hfill+D\,\Sym\,(a_1|a_2)(a_3|\omega)(a_4|\omega)
+E(a_1|\omega)(a_2|\omega)(a_3|\omega)(a_4|\omega)\Bigr)\cr
}
$$
where $D_8(c)$ is given in {\rm (\ref{eq:Dn})} and the coefficients
are listed in {\rm Appendix A.1}.
\item
{\rm (5)}
If $V$ is of class~$\CS^{10}$ then, for any $a_1,a_2,a_3,a_4,a_5\in B$, 
$$\EQ{
\Trace& \R_{a_1}\R_{a_2}\R_{a_3}\R_{a_4}\R_{a_5}\cr
&=\frac{1}{D_{10}(c)}\Bigl(
\sum A_{i_1,i_2,i_3,i_4,i_5}(a_{i_1}a_{i_2}|a_{i_3}|a_{i_4}a_{i_3})\cr
&+B_{1}\,\Cyc(a_{1}a_{2}|a_{3}a_{4})(a_5|\omega)
+B_{2}\,\Cyc(a_{1}a_{3}|a_{2}a_{4})(a_5|\omega)
+B_{3}\,\Cyc(a_{1}a_{4}|a_{2}a_{3})(a_5|\omega)\cr
&\ +C\,\Sym (a_1|a_2|a_3)(a_4|a_5)
+D\,\Sym (a_1|a_2|a_3)(a_4|\omega)(a_5|\omega)\cr
&\ +E\,\Sym (a_1|a_2)(a_3|a_4)(a_5|\omega)
+F\,\Sym (a_1|a_2)(a_3|\omega)(a_4|\omega)(a_5|\omega)\cr
&\ +G(a_1|\omega)(a_2|\omega)(a_3|\omega)(a_4|\omega)(a_5|\omega)
+H(a_1,a_2,a_3,a_4,a_5)
\Bigr)\cr
}
$$
where the summation is over the permutations of $(1,2,3,4,5)$ for
which $(a_{i_1}a_{i_2}|a_{i_3}|a_{i_4}a_{i_5})$ are distinct.  The
determinant $D_{10}(c)$ is given in {\rm (\ref{eq:Dn})} and the coefficients
are listed in {\rm Appendix A.2}.
\label{Theorem1}
\EndTheorem

\Remark
\rm
Suppose that $V$ is of class~$\CS^{8}$.
By the cyclic property of trace:
$$
\Trace R_{a_{1}}R_{a_{2}}R_{a_{3}}R_{a_{4}}
=\Trace R_{a_{2}}R_{a_{3}}R_{a_{4}}R_{a_{1}},
$$
we must have $A_{1}=A_{3}$ in Theorem \ref{Theorem1} (4) if there
exist elements $a_1,\ldots,a_4$ such that
$(a_1a_2|a_3a_4)\ne(a_1a_4|a_2a_3)$.  In this case, the dimension is
determined from the rank as
\Equation
d=-\frac{1050\ {c^6}+22565\ {c^5}+33121\ {c^4}-1707790\ {c^3}-3390408\
{c^2}+ 308160\ c}{2\ (30\ {c^5}-3212\ {c^4}+107355\ {c^3}-1135590\
{c^2}-206024\ c+825792)}.
\label{eq:constraint8}
\EndEquation
\EndRemark

We will discuss restrictions on the pair $(c,d)$ similar to
(\ref{eq:constraint8}) later in Section 3.
\subsection{Derivation by Casimir elements}
In this subsection, we sketch the derivation of the formulae in case
the form $(\ |\ )$ on $B$ is nondegenerate.

Let $\{x_1,\ldots,x_d\}$ be a basis of $B$ and let
$\{x^1,\ldots,x^d\}$ be the dual basis with respect to the form $(\ |\
)$: $(x_i|x^j)=\delta_{i,j}$.  We suppose that any expression with a
repeated index $i$ must be summed over $i=1,\ldots,d$.
Note that
\Equation
(x^i_{(3)}a)_{(-1)}x_i=x^i_{(-1)}a_{(3)}x_i=a.
\label{eq:ONB}
\EndEquation

The strategy of the derivation of the formula is to write the trace as
$$
\Trace \R_{a_1}\cdots\R_{a_m}
=(a_1{}_{(1)}\cdots a_{m}{}_{(1)}x_i|x^i)
=(x^i_{(1)}a_{m-1}{}_{(1)}\cdots a_1{}_{(1)}x_i|a_m),
$$
and to rearrange the vectors by using the identity
(\ref{eq:borcherds}) and the invariance of the form $(\ |\ )$ on $V$
until the trace is written as a sum of expressions of the form
$(x^i_{(k)}x_i|X)$ where $X$ is an element of $V$ written by
$a_1,a_2,\ldots,a_m$ and $k$ is an integer.

Let us first study the `Casimir' elements:
\Equation
\casimir_{n}=x^i_{(3-n)}x_i=\sum_{i=1}^d x^i_{(3-n)}x_i,
\EndEquation
which do not depend on the choice of the basis.
Note that
\Equation
\casimir_0=d\1\hbox{\quad and \quad} \casimir_1=0.
\label{eq:init}\EndEquation
and that the vector $\casimir_{n}$ for an odd $n$ is determined from
those for even $n$ by the action of $L_{-1}$.  By the identity
(\ref{eq:borcherds}), the sequence $\casimir_0,\ldots,\casimir_n$ is
subject to the relations
\Equation
L_m\casimir_k=(m+k-2)\casimir_{k-m}+\delta_{m,2}L_{-k+2}\1+\delta_{m,k-2}\frac{m^3-m}{6}L_{-2}\1,\ (k=0,\ldots,n).
\label{eq:cond}\EndEquation

Now, in order to deduce some information on the Griess algebra by
these elements, the following simple observation is fundamental.
\Lemma
\sl
The vector $\casimir_{n}=x^i_{(3-n)}x_i$ is fixed by any automorphism of the VOA $V$.
\EndLemma
%\Proof
%
%For any automorphism $g$ of the VOA $V$, we have
%$(g(x_i)|g(x^j))=(x_i|x^j)=\delta_{ij}$. Hence the set
%$\{g(x_1),\ldots,g(x_d)\}$ is a basis and $\{g(x^1),\ldots,g(x^d)\}$
%is the dual basis.  Thus $g(x^i_{(3-n)} x_i)=
%g(x^i)_{(3-n)}g(x_i)=x^i_{(3-n)} x_i$.
%
%\EndProof

Therefore, if the VOA $V$ has larger symmetry then the vector
$\casimir_n$ has to belong to a smaller subspace of $V$; the smallest
possible case is $V_\omega$ when the VOA is of class~$\CS^n$.

\Lemma
\sl
If the VOA $V$ is of class~$\CS^n$ then the Casimir elements
$\casimir_2,\ldots,\casimir_n$ are contained in $V_\omega$.
\EndLemma

Now, suppose that the vectors $\casimir_2,\ldots,\casimir_n$ are
indeed contained in $V_\omega$.  If the Virasoro submodule $V_\omega$
does not contain a singular vector of degree up to $n$, then these
vectors are uniquely determined by the properties (\ref{eq:init}) and
(\ref{eq:cond}).  In this way, we have the following result.

\Prop
\sl
If the VOA $V$ is of class~$\CS^n$ then $\casimir_n$ is uniquely written as
\Equation
\casimir_n=\frac{1}{D_n(c)}\sum_{\vec{m}\in P_n} P_{\vec{m}}(c,d)\basis{\vec{m}}.
\label{EqSolution}\EndEquation
where $P_{\vec{m}}(c,d)$ are certain polynomials in $c$ and $d$.
\EndProp

The explicit expressions of $\casimir_n$ for $n=2,4,6$ are given as follows:
\Equation
\EQ{
\casimir_2=&\frac{4d}{c}\basis{2},\quad\casimir_4=\frac{6(d-1)}{5c+22}\basis{4}+\frac{2(5c+22d)}{c(5c+22)}\basis{2,2},\cr
\casimir_6=&\frac{8(d-1)(5c^2+35c-228)}{(2c-1)(5c+22)(7c+68)}\basis{6}\cr
&+\frac{2c(70c^2+769c+1644)+4d(92c^2+427c-748)}{c(2c-1)(5c+22)(7c+68)}
\basis{4,2}\hfill\cr
&+\frac{31(d-1)(5c+44)}{(2c-1)(5c+22)(7c+68)}\basis{3,3}
+\frac{992(d-1)}{(2c-1)(5c+22)(7c+68)}\basis{2,2,2}.\hskip-3em\hfill\cr
}
\label{eq:explicit}
\EndEquation
The expressions for higher $n$ are so lengthy; we do not include them
in this paper.  The expressions in case $c=24$ and $d=196884$ look as follows:
{\small
$$
\EQ{
\casimir_2&=32814\basis{2},\quad \casimir_4=8319\basis{4}+2542\basis{2,2},\cr
\casimir_6&=3492\basis{6}+1302\basis{4,2}+\frac{1271}{2}\basis{3,3}+124\basis{2,2,2},\cr
\casimir_8&=\frac{3863}{2}\basis{8}+552 \basis{6,2}+434\basis{5,3}+\frac{333}{2}\basis{4,4}+96\basis{4,2,2}+93\basis{3,3,2}+\frac{13}{3}\basis{2,2,2,2},\cr
\casimir_{10}&=1182\basis{10}+\frac{613}{2}\basis{8,2}+207\basis{7,3}+141\basis{6,4}+41\basis{6,2,2}+74\basis{5,5}+64\basis{5,3,2}\cr
&+\frac{99}{4}\basis{4,4,2}+24\basis{4,3,3}+\frac{9}{2}\basis{4,2,2,2}+\frac{13}{2}\basis{3,3,2,2}+\frac{7}{60}\basis{2,2,2,2,2}.\cr
}
$$
}

Now, suppose that $V$ is of class~$\CS^2$ and let $a$ be any element of the
Griess algebra.  Then we immediately obtain Theorem \ref{Theorem1} (1):
$
\Trace \R_{a}
=(a_{(1)}x^i|x_i)
=(x^i_{(1)}x_i|a)
={4d}(a|\omega)/c.
$
In particular, we have
\Equation
\Trace R_{ab}=\frac{8d}{c}(a|b).
\EndEquation

Next, take any two elements $a,b$ of the Griess algebra.  Then
$$
\Trace \R_{a}\R_{b}
=(a_{(1)}b_{(1)}x^i|x_i)
=(x^i_{(1)}a_{(1)}x_i|b).
$$
By the identity (\ref{eq:borcherds}) for $p=-1,q=1,r=2$, we have
$
-(x^i_{(3)}a)_{(-1)}x_i
= x^i_{(1)}a_{(1)}x_i+x^i_{(-1)}a_{(3)}x_i
-a_{(3)}x^i_{(-1)}x_i+2a_{(2)}x^i_{(0)}x_i-a_{(1)}x^i_{(1)}x_i.
$
Since $a_{(2)}x^i_{(0)}x_i=a_{(1)}x^i_{(1)}x_i$, 
$$\EQ{
\Trace \R_{a}\R_{b}
&=-2(a|b)+(a_{(3)}x^i_{(-1)}x_i|b)-(a_{(1)}x^i_{(1)}x_i|b)\cr
&=-2(a|b)+(x^i_{(-1)}x_i|a_{(-1)}b)-(x^i_{(1)}x_i|a_{(1)}b)\cr
}
$$
Suppose that $V$ is of class~$\CS^4$ and substitute the expressions of
$\kappa_2=x^i_{(1)}x_i$ and $\kappa_4=x^i_{(-1)}x_i$ given by
(\ref{eq:explicit}). Then using
\Equation
\EQ{
&(L_{-4}\1|a_{(-1)}b) =(\1|L_4a_{(-1)}b) =6(a|b),\cr
&(L_{-2}L_{-2}\1|b_{(-1)}a)
=(\1|L_2L_2(b_{(-1)}a))
=2(a|\omega)(b|\omega)+8(a|b),\cr
}
\EndEquation
we have Theorem \ref{Theorem1} (2):
$$
\EQ{
&\Trace \R_{a}\R_{b}\cr
&\quad=-2(a|b)+\frac{6(d-1)}{5c+22}6(a|b)+\frac{2(5c+22d)}{c(5c+22)}\Bigl(2(a|\omega)(b|\omega)+8(a|b)\Bigr)-\frac{4d}{c}(ab|\omega)\cr
&\quad=\frac{-2(5c^2-88d+2cd)}{c(5c+22)}(a|b)+\frac{4(5c+22d)}{c(5c+22)}(a|\omega)(b|\omega)
}$$

Similarly, one can obtain a expression in terms of the inner product
and the multiplication of the trace in which three and four elements
of the Griess algebra are concerned if $V$ is of class~$\CS^6$ and
of~$\CS^8$, respectively.  However, if five elements are concerned,
then we encounter the expression
$
a_1{}_{(3)}a_2{}_{(2)}a_3{}_{(1)}a_4{}_{(0)}a_5
$ 
and its permutations which cannot be written by a combination of the
inner product $(\ |\ )$ and the multiplication in general.  Thus we
are led to consider the totally antisymmetric quinery form defined by
(\ref{eq:quinery}).  Then the trace $\Trace
\R_{a_1}\R_{a_2}\R_{a_3}\R_{a_4} \R_{a_5}$ is written by a
combination of these operations if $V$ is of class~$\CS^{10}$.
In this way, we obtain Theorem \ref{Theorem1} (3)--(5).
\subsection{Derivation by projection to $V_\omega$}
In this subsection, we will sketch another derivation of the formulae
under the assumption that $\Aut V$ is finite.

For any $v\in V$, let $\tildeb{v}$ denote its average over the action
of the automorphism group:
\Equation
\tildeb{v}=\frac{1}{|\Aut V|}\sum_{g\in \Aut V}gv.
\EndEquation
The following lemma is obvious.
\Lemma
\sl
Let $v$ be an element of $V^n$. Then 
$
\Trace|_{V^k} v_{(n-1)}=\Trace|_{V^k} (gv)_{(n-1)}
$
for any $g\in\Aut V$ at any degree $k$.
\label{lemma2.6}
\EndLemma
In particular, we have $\Trace v_{(n-1)}=\Trace \tildeb{v}_{(n-1)}$
for any $v\in V^n$.

Now, for each $n$, consider the map
\Equation
\eta_n:V^n\rightarrow \CC^{\,\dim V^n},\quad v\mapsto (L_{\vec{m}}v)^{}_{\vec{m}\in P_n},
\EndEquation
where $P_n$ is the set (\ref{eq:partition}) which parametrizes a basis
of $V_\omega^n$.

\Lemma
\sl
If $D_n(c)\ne0$ then
$
V^n=V_\omega^n\oplus \Ker\eta_n.
$
\EndLemma
\Proof Since the map $\eta_n$ is isomorphic on $V_\omega^n$ if
$D_n(c)\ne0$, the result follows.
\EndProof

For any $v\in V$, let $\delta(v)$ denote its projection to $V_\omega$ with respect to the decomposition as in the lemma.

\Lemma
\sl
Suppose that $V$ is of class~$\CS^n$.
Then
$
\Trace|_{V^k} v_{(n-1)}=\Trace|_{V^k} \delta(v)_{(n-1)}
$
for any $v\in V^n$ at any degree $k$.
\label{lemma:delta}
\EndLemma
\Proof
Set $\delta=\delta(v)$ and $\pi=v-\delta(v)$.  Then obviously
$
\tildeb{v}=\tildeb{\delta}+\tildeb{\pi}=\delta+\tildeb{\pi},
$
and $\tildeb{\pi}$ is fixed by any automorphism of $V$.  Therefore, if
$V$ is of class~$\CS^n$ then $\tildeb{\pi}=0$ because $\tildeb{\pi}\in
V_\omega\cap \Ker\eta_n$.  Hence, by Lemma \ref{lemma2.6}, we have
$
\Trace|_{V^k} v_{(n-1)}
=\Trace|_{V^k} \tildeb{v}_{(n-1)}
=\Trace|_{V^k} \delta_{(n-1)}.
$
\EndProof

Since $\delta(v)\in V_\omega$, the action $\delta(v)_{(n-1)}$ on $B$
can be explicitly computed by the commutation relation of the Virasoro
algebra.

For example, if $c\ne 0$, then we have
$
\delta(a)=2(a|\omega)\omega/c
$
for any $a\in B$.  Namely, any element $a\in B$ can be written as
\Equation
a=\frac{2(a|\omega)}{c}\omega+\pi,
\label{eq:a=}
\EndEquation
where $\pi$ is a primary vector in $B=V^2$.  If $V$ is of class~$\CS^2$
then 
$$
\Trace \R_a
=\frac{2(a|\omega)}{c}\Trace \R_\omega
=\frac{4d}{c}(a|\omega).
$$
Thus we have obtained Theorem \ref{Theorem1} (1).

Next, let $a,b$ be any elements of $B$. 
Then, by the identity (\ref{eq:borcherds}) for $p=2$, $q=1$, $r=-1$, we have 
$
a_{(1)}b_{(1)}=(a_{(-1)}b)_{(3)}+2(a_{(0)}b)_{(2)}+(a_{(1)}b)_{(1)}-a_{(-1)}b_{(3)}-b_{(-1)}a_{(3)}
$
on $B$.
Therefore, 
$$
\Trace \R_a\R_b=\Trace|_B (a_{(-1)}b)_{(3)}+2\Trace|_B (a_{(0)}b)_{(2)}+\Trace|_B(a_{(1)}b)_{(1)}-2(a|b).
$$
If $V$ is of class~$\CS^4$ then, since
\Equation
\EQ{
\hfill\delta(a_{(-1)}b)
&=&\frac{6c(a|b)-12(a|\omega)(b|\omega)}{c(5c+22)}\basis{4}
+\frac{44(a|b)+20(a|\omega)(b|\omega)}{c(5c+22)}\basis{2,2},\cr
\hfill\delta(a_{(0)}b)
&=&\frac{2(a|b)}{c}\basis{3},\quad\delta(a_{(1)}b)=\frac{4(a|b)}{c}\basis{2}\cr
}
\label{eq:projection}
\EndEquation
and 
$\Trace|_B [4]_{(3)}=6d$, $\Trace|_B [2,2]_{(3)}=8d+c$, $\Trace|_B
[3]_{(2)}=-4d$, $\Trace|_B [2]_{(1)}=2d$,
we have Theorem \ref{Theorem1} (2).

The derivation of Theorem \ref{Theorem1} (3)--(5) is similar.
\section{Constraints on $c$ and $d$}
Let $B$ be the Griess algebra of a VOA $V$ such that the form $(\ |\
)$ on $B$ is nondegenerate.

In this section, we will give some necessary conditions satisfied by
the pair $(c,d)$ for a VOA with larger symmetry under some additional
conditions on $V$.
\subsection{Constraints from a proper idempotent}
By a {\it proper idempotent}\/ we mean an idempotent $e\in B$ such
that the central charge $b=2(e|e)$ differs from $0$ and $c$.  If the
algebra $B$ has a real form on which the bilinear form $(\ |\ )$ is
positive-definite then, by Theorem 11 of \cite{MeyerNeutsch} and
Theorem 6.8 of \cite{Miyamoto1}, the conformal vector $\omega$ is
decomposable if $d\geq2$; In particular, the algebra $B$ contains a
proper idempotent.

\Lemma
\sl
Let $\varphi(\idem)$ be a vector generated by an idempotent $\idem\in B$.  Then
$$
(\idem_{(2)}(a_{(m-1)}b)|\varphi(\idem))=(3-m)(a_{(m)}b|\varphi(\idem))
$$
holds for any $a,b\in B$ and $m\in\ZZ$.
\label{Lemma3}
\EndLemma
\Proof
Note that the actions of $\idem_{(p)}$ and $(\omega-\idem)_{(q)}$ commute, and
$(\omega-\idem)_{(q)}\1=0$ if $q\geq 0$.
Hence $\idem_{(q)}\varphi(\idem)=\omega_{(q)}\varphi(\idem)$ if $q\geq 0$.
Therefore, by the invariance of the bilinear form $(\ |\ )$ on $V$, 
$
(\idem_{(2)}(a_{(n)}b)|\varphi(\idem))
=(a_{(n)}b|\idem_{(0)}\varphi(\idem))
=(a_{(n)}b|\omega_{(0)}\varphi(\idem))
=(\omega_{(2)}(a_{(n)}b)|\varphi(\idem)).
$
The result follows from
$
(\omega_{(2)}(a_{(m-1)}b)|\varphi(\idem))
=((\omega_{(0)}a)_{(m+1)}b)|\varphi(\idem))
+2((\omega_{(1)}a)_{(m+1)}b)|\varphi(\idem))
=(3-m)(a_{(m)}b)|\varphi(\idem)).
$
\EndProof

Now suppose that $V$ is of class~$\CS^6$.  By the lemma, 
\Equation
(x^i_{(-1)}\idem_{(1)}x_i|\idem_{(0)}\idem_{(0)}\idem)=3(x^i_{(0)}\idem_{(1)}x_i|\idem_{(0)}\idem).
\EndEquation
Computing the both-hand sides of this equality by the method of
subsection 2.2, we have
$$
b(b-c)\Bigl((70{c^2}+955{c}+2388)c-2({c^2}-55c+748)d\Bigr)=0,
$$
where $b=2(\idem|\idem)$ is the central charge of $\idem$.  Therefore,
if $\idem$ is proper then
\Equation 
d =\frac{(70{c^2}+955{c}+2388)c}{2({c^2}-55c+748)}
=35c+2402+\frac{c^2+214248c-3593392}{2(c^2-55c+748)}.
\label{eq:constraint6}
\EndEquation

Further, if $V$ is of class~$\CS^8$ then, by computing
\Equation
(x^i_{(-3)}e_{(1)}x_i|e_{(0)}e_{(0)}e_{(0)}e_{(0)}e)
=5(x^i_{(-2)}e_{(1)}x_i|e_{(0)}e_{(0)}e_{(0)}e),
\EndEquation
we have
\Equation
d=\frac{5250{c^5}+155250{c^4}+1369715{c^3}+3507098{c^2}+1497768c}
    {125{c^4}-4770{c^3}-23382{c^2}+1561868c+1032240}
\label{eq:constraint8'}
\EndEquation
if $\idem$ is proper.

\Theorem \sl 
Let $V$ be a VOA of class~$\CS^8$ for which the form $(\ |\ )$ on $B$
is nondegenerate.  If $B$ contains a proper idempotent then $c=24$ and
$d=196884$.
\label{Theorem2}
\EndTheorem
\Proof
By (\ref{eq:constraint6}) and (\ref{eq:constraint8'}),
$
c=-46/3,\ -68/7,\ -22/5,\ -3/5,\ 0,\ 1/2,\ 24,\ 142/5.
$
However, $D_6(c)=0$ for the first 6 cases and $d<0$ for the last case.
\EndProof

By the remark at the beginning of this subsection, we have the following corollary.

\Cor\sl 
Suppose that the algebra $B$ has a real form on which the bilinear
form $(\ |\ )$ is positive-definite.  If the VOA $V$ is of
class~$\CS^8$ and if $d\geq 2$ then $c=24$ and $d=196884$.
\EndCor

Now, let us come back to VOA's of class~$\CS^6$ but restrict our
attention to the case when the rank $c$ is a positive half-integer.
In this case, inspecting the relation (\ref{eq:constraint6}), we see
that the pair $(c,d)$ must be one from the following table:
$$
\EQ{
\EQ{
&\hfill c&\quad&\hfill\hfill\hfill d\hfill
&\qquad&\hfill c&\quad&\hfill\hfill d\hfill
&\qquad&\hfill c\hfill&\quad&\hfill\ d\hfill
&\qquad&\hfill c&\quad&\hfill\  d \hfill\cr
\noalign{\vskip1ex}
&\hfill 8&&\hfill 156&&\hfill 23&\onehalf&\hfill 96256&&\hfill 32&&\hfill 139504&&\hfill 54&\onehalf&\hfill 9919\cr
&\hfill 16&&\hfill 2296&&\hfill 24&&\hfill 196884&&\hfill 34&&\hfill 57889&&\hfill 68&&\hfill 8146\cr
&\hfill 20&&\hfill 10310&&\hfill 24&\onehalf&\hfill 1107449&&\hfill 36&&\hfill 35856&&\hfill 93&\onehalf&\hfill 7566\cr
&\hfill 21&\onehalf&\hfill 21414&&\hfill 30&\onehalf&\hfill 1964871&&\hfill 40&&\hfill 20620&&\hfill 132&&\hfill 8154\cr
&\hfill 22&&\hfill 28639&&\hfill 31&\onehalf&\hfill 207144&&\hfill 44&&\hfill 14994&&\hfill 1496&&\hfill 54836\cr
}\cr
\noalign{\vskip3ex}
\hbox{\bf Table 3.1}\cr
\noalign{\vskip2ex}
}
$$
\subsection{Constraints from an idempotent of central charge $1/2$}
Suppose that $B$ contains an idempotent of central charge $1/2$ for
which the eigenvalues of the adjoint action are $0,1/2,1/16$ and $2$,
and the eigenspace with eigenvalue $2$ is one-dimensional.  This is
indeed the case if $e$ generates a subVOA isomorphic to
$\Lvir(\frac12,0)$, for this is a rational VOA for which the
irreducible modules are isomorphic to either $\Lvir(\frac12,0)$,
$\Lvir(\frac12,\frac12)$ or $\Lvir(\frac12,\frac1{16})$.  In
particular, this holds if $V$ has a real form on which the form $(\ |\
)$ is positive-definite and $e$ is a real idempotent of central charge
$1/2$.

First consider the case when the eigenspace with eigenvalue $1/16$ is
zero.  In this case, we have the following equations satisfied by
the dimension $d(\frac12)$ of the eigenspace:
\Equation
d(\textstyle\frac12)+4=2\Trace \R_{e},\ d(\textstyle\frac12)+16=4\Trace \R_{e}^2
\EndEquation
By the compatibility, we obtain
$
2\Trace \R_{e}^2-\Trace \R_{e}=6.
$
If the VOA $V$ is of class~$\CS^4$ then, by Theorem \ref{Theorem1}, 
\Equation
(-22 + 2 c) d = (-37 - 10 c).
\label{eq:constraint}
\EndEquation
Further if $V$ is of class~$\CS^6$ then we get
\Equation
(3{c^2}+164c-2992)d=c(-140{c^2}-1903c-4832).
\EndEquation
By solving these equations, we obtain the following result.

\Prop
\sl
Suppose that the VOA $V$ contains an idempotent of central charge
$1/2$ for which the eigenvalues are $0,1/2$ and $2$ and the eigenspace
with eigenvalue~$2$ is one-dimensional. If $V$ is of class~$\CS^6$ and
$d\geq 2$ then $c=8$ and $d=156$.  
\EndProp

Now suppose that the rank $c$ is a positive half-integer.  If $V$ is
of class~$\CS^4$ then, by (\ref{eq:constraint}), the nonnegativity and
the integrality of $d(0)$ and $d(\frac12)$ gives us a list of possible
pairs of such $c$ and $d$.  They are given by the following table:
$$
\EQ{
\EQ{
&\hfill c&&\hfill d&\;=\;&\hfill d(0)\;+\;&\hfill d(\textstyle\frac12)\;+\;1\cr
\noalign{\vskip1ex}
&\hfill 0&\onehalf\quad &\hfill 1&\;=\;&\hfill 0\;+\;&\hfill 0\;+\;1\cr
&\hfill 4&&\hfill22 &\;=\;&\hfill 14\;+\;&\hfill 7\;+\;1\cr
&\hfill 7&\onehalf\quad &\hfill120&\;=\;&\hfill 91\;+\;&\hfill 28\;+\;1\cr
&\hfill 8&&\hfill 156&\;=\;&\hfill 120\;+\;&\hfill 35\;+\;1\cr
&\hfill 9&\onehalf &\hfill418&\;=\;&\hfill 333\;+\;&\hfill 84\;+\;1\cr
&\hfill 10&&\hfill 685&\;=\;&\hfill 551\;+\;&\hfill 133\;+\;1\cr
&\hfill 10&\onehalf&\hfill 1491&\;=\;&\hfill 1210\;+\;&\hfill 280\;+\;1\cr
}\cr
\noalign{\vskip3ex}
\hbox{\bf Table 3.2}\cr
\noalign{\vskip2ex}
}
$$

\Note\rm
The fixed-point VOA $V_{\sqrt{2}E_8}^{\scriptscriptstyle +}$ of the
VOA associated with the lattice $\sqrt{2}E_8$ by the $-1$ automorphism
of the lattice would be an example with $c=8$ and $d=156$.  According
to \cite{Griess2}, the automorphism group of this VOA is isomorphic to
$O^+_{10}(2)$.  The decomposition $156=120+35+1$ coincides with
Theorem 5.2 of \cite{Griess2}.  The Hamming code VOA $V_{H_8}$
considered by Miyamoto \cite{Miyamoto2}, isomorphic to the VOA
$V^+_{\sqrt{2}D_4}$, would be an example with $c=4$ and $d=22$. The
automorphism group of this VOA is isomorphic to a group of shape
$2^6:(\GL_3(2)\times S_3)$ \cite{MatsuoMatsuo}.
\EndNote

We next consider the case when the $1/16$ components are indeed present.
In this case, if $V$ is of class~$\CS^6$ then we have
\Equation
(2c^2-110c+1496)d=(70c^2+955c+2388)c,
\EndEquation
which is actually the same as the condition (\ref{eq:constraint6}).

If $V$ is of class~$\CS^6$ with $d\geq 2$ and if $c$ is a positive
half-integer then the rank $c$ and the dimension $d$ must be a pair
from the following table:
$$
\EQ{
\EQ{
&\hfill c&&\hfill d&\;=\;\hfill d(0)\;+\;&\hfill d(\textstyle\frac12)\;+\;&\hfill d(\textstyle\frac1{16})\;+\;1\cr
\noalign{\vskip1ex}
&\hfill 16&&\hfill  2296&\;=\;\hfill 1116\;+\;&\hfill 155\;+\;&\hfill 1024\;+\;1\cr
&\hfill 20&&\hfill  10310&\;=\;\hfill 4914\;+\;&\hfill 403\;+\;&\hfill 4992\;+\;1\cr
&\hfill 23&\onehalf\quad&\hfill  96256&\;=\;\hfill 46851\;+\;&\hfill 2300\;+\;&\hfill 47104\;+\;1\cr
&\hfill 24&&\hfill  196884&\;=\;\hfill 96256\;+\;&\hfill 4371\;+\;&\hfill 96256\;+\;1\cr
&\hfill 24&\onehalf&\hfill  1107449&\;=\;\hfill 543960\;+\;&\hfill 22816\;+\;&\hfill 540672\;+\;1\cr
&\hfill 30&\onehalf&\hfill  1964871&\;=\;\hfill 1029630\;+\;&\hfill 13640\;+\;&\hfill 921600\;+\;1\cr
&\hfill 31&\onehalf&\hfill  207144&\;=\;\hfill 109771\;+\;&\hfill 1116\;+\;&\hfill 96256\;+\;1\cr
&\hfill 32&&\hfill  139504&\;=\;\hfill 74340\;+\;&\hfill 651\;+\;&\hfill 64512\;+\;1\cr
&\hfill 36&&\hfill  35856&\;=\;\hfill 19951\;+\;&\hfill 0\;+\;&\hfill 15904\;+\;1\cr
}\cr
\noalign{\vskip3ex}
\hbox{\bf Table 3.3}\cr
\noalign{\vskip2ex}
}$$

\Note\rm 
The moonshine module $V^\natural$ is of course an example with $c=24$
and $d=196884$.  The bosonic projection of the Babymonster VOSA
$V\!B^\natural$ constructed by H\"ohn \cite{Hohn} would be an example
with $c=23\onehalf$ and $d=96256$.  The fixed-point VOA
$V^{\scriptscriptstyle +}_{\Lambda_{16}}$ of the VOA associated with
the Barnes-Wall lattice $\Lambda_{16}$ by the $-1$ automorphism of the
lattice would be an example with $c=16$ and $d=2296$.
\EndNote
\section{Application to the Moonshine Module}
Now, let $V^\natural$ be the moonshine module and let $B^\natural$ be
the Conway-Griess algebra.  In this section, we will compute the
spectrum of the eigenspace decomposition of $B^\natural$ with respect
to idempotents related to some Monster elements starting from the
trace formulae by the representation theory of various subVOA's inside
$V^\natural$.
\subsection{Norton's formula}
Since the moonshine module $V^\natural$ is of class~$\CS^{11}$, just
substituting $c=24$ and $d=196884$ in Theorem \ref{Theorem1}, we
recover the original trace formulae of Norton \cite{Norton}:
\Cor\sl
For any elements $a_1,a_2,a_3,a_4,a_5$ of the Conway-Griess algebra
$B^\natural$, 
{
$$
\EQ{
\hfill\Trace \R_{a_1}&=32814 (a_1|\omega),\cr
\hfill\Trace \R_{a_1}\R_{a_2}&=4620(a_1|a_2)+5084(a_1|\omega)(a_2|\omega),\cr
\hfill\Trace \R_{a_1}\R_{a_2}\R_{a_3}&=900 (a_1|a_2|a_3)
+620 \,\Cyc\,(a_1|a_2)(a_3|\omega)
+744 (a_1|\omega)(a_2|\omega)(a_3|\omega),\cr
\hfill\Trace \R_{a_1}\R_{a_2}\R_{a_3}\R_{a_4}&=
166(a_1 a_2|a_3 a_4)-116(a_1 a_3|a_2 a_4)+166(a_1 a_4|a_2 a_3)\cr
&+114 \,\Sym\, (a_1|a_2|a_3)(a_4|\omega)
+52 \,\Sym\, (a_1|a_2)(a_3|a_4)\cr
&+80 \,\Sym\,(a_1|a_2)(a_3|\omega)(a_4|\omega)
+104 (a_1|\omega)(a_2|\omega)(a_3|\omega)(a_4|\omega),\cr
\hfill\Trace \R_{a_1}\R_{a_2}\R_{a_3}\R_{a_4}\R_{a_5}
&=30\Cyc\,(a_1 a_2|a_3|a_4 a_5)
+4\Cyc\,(a_1 a_4|a_3|a_2 a_5)
-22\Cyc\,(a_1 a_5|a_3|a_2 a_4)\hskip-1ex\cr
&\hskip-10ex+20\Cyc\,(a_1 a_2|a_3 a_4)(a_5|\omega)
-14\Cyc\,(a_1 a_3|a_2 a_4)(a_5|\omega)
+20\Cyc\,(a_1 a_4|a_2 a_3)(a_5|\omega)\hskip-1ex\cr
&\hskip-5ex+\,8\,\Sym\,(a_1|a_2|a_3)(a_4|a_5)
+14\,\Sym\,(a_1|a_2|a_3)(a_4|\omega)(a_5|\omega)\cr
&\hskip-5ex+\,6\,\Sym\,(a_1|a_2)(a_3|a_4)(a_5|\omega)
+10\,\Sym\,(a_1|a_2)(a_3|\omega)(a_4|\omega)(a_5|\omega)\cr
&\hskip-5ex+14(a_1|\omega)(a_2|\omega)(a_3|\omega)(a_4|\omega)(a_5|\omega)+52(a_1,a_2,a_3,a_4,a_5).\cr}
$$
}\label{Corollary}
\EndCor
\Note\rm
To compare this result with Table 2 of \cite{Norton}, substitute
$1=\omega/2$ (the identity element of the algebra) and suppose for the
last one that $a_1,a_2,a_3,a_4,a_5$ are perpendicular to $\omega$.
The formula for $\Trace \R_{a_1}\R_{a_2}$ was given in p.\ 528
of \cite{Conway}.
\EndNote

In particular, letting $a_1=\cdots=a_5$ be an idempotent $\idem$, we
have
\Equation
\EQ{
\Trace \,\R_{\idem}\hfill&=&\hfill  32814\,(\idem|\idem)&,\hfill\cr
\Trace \R_{\idem}^{2}\hfill&=&\hfill 4620\,(\idem|\idem)&+&\hfill 5084\,{(\idem|\idem)^2}&,\hfill\cr
\Trace \R_{\idem}^{3}\hfill&=&\hfill 1800\,(\idem|\idem)&+&\hfill 1860\,{(\idem|\idem)^2}&+&\hfill 744\,{(\idem|\idem)^3},&\hfill\cr
\Trace \R_{\idem}^{4}\hfill&=&\hfill 864\,(\idem|\idem)&+&\hfill 1068\,{(\idem|\idem)^2}&+&\hfill 480\,{(\idem|\idem)^3}&+&\hfill 104\,{(\idem|\idem)^4}&,\hfill\cr
\Trace \R_{\idem}^{5}\hfill&=&\hfill 480\,(\idem|\idem)&+&\hfill 680\,{(\idem|\idem)^2}&+&\hfill 370\,{(\idem|\idem)^3}&+&\hfill 100\,{(\idem|\idem)^4}&+&\hfill 14\,{(\idem|\idem)^5}.\cr
}
\label{Cor}
\EndEquation
%$$
%\EQ{
%&\Trace \R_{\idem}&=16407 b,\cr
%&\Trace \R_{\idem}^2&=1271b^2+2310 b,\cr
%&\Trace \R_{\idem}^3&=93b^3+465 b^2 +900 b,\cr
%&\Trace \R_{\idem}^4&=\frac{13}{2}b^4+60 b^3 +267 b^2+432b,\cr
%&\Trace \R_{\idem}^5&=\frac{7}{16}b^5+\frac{25}{4}b^4+\frac{185}{4}b^3+170b^2+240b,\cr
%}
%$$
Recall that $2\,(\idem|\idem)$ is the central charge of the Virasoro algebra corresponding to $\idem$.

\Remark \rm 
By a slight variation of the argument of Subsection 2.3, using 
results of Zhu \cite{Zhu} and Dong and Mason \cite{DongMason} and the
absence of a cusp form of weight less than $12$, we see that Norton's
formulae hold for any rational selfdual (holomorphic) VOA of rank $24$
with shape (\ref{EqDoublyBosonic}) satisfying the $C_2$ finiteness
condition.
\EndRemark
\subsection{Eigenspace decomposition of the Conway-Griess algebra}
Let $U\rightarrow V^\natural$ be an inclusion of a VOA into
$V^\natural$ such that
\vskip0.5ex\noindent
(1) The map preserves the operations of VOA's.
\vskip0.5ex\noindent
(2) The vacuum vector of $U$ is mapped to the vacuum vector of $V^\natural$.
\vskip0.5ex\noindent
(3) The conformal vector $t$ of $U$ is mapped to an idempotent $\idem$
in $B^\natural$.
\vskip0.5ex\noindent
(4) The action $\idem_{(n)}$ for nonnegative $n$ coincides with that
of $\omega_{(n)}$ on the image of $U$.
\vskip0.5ex\noindent

Recall that $V^\natural$ has a real form $V^\natural_\RR$ on which the
form $(\ |\ )$ is positive-definite.  We suppose that $\idem$ above is
a real idempotent.  Then the adjoint action $R_\idem=\idem_{(1)}$ is
semisimple on $B^\natural$. Let $B^\natural(\lambda)$ denote the
eigenspace of $\R_\idem$ with eigenvalue $\lambda$, and let
$d(\lambda)$ be its dimension.

For the representation theory of various VOA's discussed below, we
refer the reader to \cite{FrenkelZhu}, \cite{Wang}, \cite{DMZ},
\cite{KMY}, \cite{Miyamoto4}, \cite{DongNagatomo}, \cite{Abe} as well
as the physics papers \cite{BPZ}, \cite{FateevZamolodchikov},
\cite{ZamolodchikovFateev}, \cite{FateevLukyanov}, and for the
construction of an automorphism by means of fusion rules to
\cite{Miyamoto1} and \cite{Miyamoto4}.  We will use the ATLAS notation
\cite{ATLAS} for conjugacy classes of the Monster.
\subsubsection*{$\Lvir(\frac12,0)$ and 2A involution}
Let $U$ be the VOA $\Lvir(\frac12,0)$ associated with the irreducible
highest weight representation of the Virasoro algebra at $c=1/2$.
There are $3$ irreducible modules for this VOA, which are parametrized
by the lowest conformal weight $h=0,\ 1/16,\ 1/2$.  Hence we have a
decomposition
$$
\textstyle B^\natural=B^\natural(0)\oplus B^\natural(\frac1{16})\oplus
B^\natural(\frac12)\oplus B^\natural(2),
$$
where $d(2)=1$. Hence
$$\frac{d({\textstyle\frac12})}{2}+\frac{d({\textstyle\frac1{16}})}{16}+2=\Trace \R_{\idem},
\quad 
\frac{d({\textstyle\frac12})}{2^2}+\frac{d({\textstyle\frac1{16}})}{16^2}+2^2=\Trace \R_{\idem}^2.
$$
By (\ref{Cor}) with $2(\idem|\idem)=1/2$, we get
$$
d(0)=d({\textstyle\frac{1}{16}})=96256,\
d({\textstyle\frac{1}{2}})=4371,\ d(2)=1.
$$

Now consider the map
$$
\EQ{
\hfill 1
&\quad\hbox{on}\quad &B^\natural(0)\oplus B^\natural({\textstyle\frac{1}{2}})\oplus B^\natural(2),
\quad
-1
&\quad\hbox{on}\quad &B^\natural({\textstyle\frac{1}{16}}).
\cr
}
$$
This map gives rise to an automorphism of $B^\natural$
\cite{Miyamoto1}, so an involution of the Monster.  It is identified
with a 2A involution of the Monster by \cite{Miyamoto1} and
\cite{Conway}.  We may confirm this by looking at the trace of this
map; it is $96256-96256+4371+1=4372$, which coincides with the
corresponding value in the list of Conway and Norton
\cite{ConwayNorton}.

Thus we have come back to the situation considered in Section 15 of
\cite{Conway} without using any explicit construction of the
Conway-Griess algebra.
\subsubsection*{$\Lvir(\frac12,0)\otimes \Lvir(\frac12,0)$ and 2B involution}
Suppose given an embedding of $U=\Lvir(\frac12,0)\otimes
\Lvir(\frac12,0)$ into $V^\natural$, and let $e^1$ and $e^2$ denote
the images of the conformal vector of the first and the second
component which we suppose to be real. Since they are mutually
orthogonal, we have a simultaneous eigenspace decomposition
$$
B^\natural=\bigoplus_{h,h'\in\{0,\frac1{16},\frac12,2\}}
B^\natural(h,h').
$$
We already know that $d(2,0)=d(0,2)=1$ and
$d(2,\frac1{16})=d(2,\frac12)=d(\frac1{16},2)=d(\frac12,2)=d(2,2)=0$.
By Corollary \ref{Corollary}, we have
$$
\Trace \R_{e^1}\R_{e^2}=\frac{1271}{4},\ 
\Trace \R_{e^1}^2\R_{e^2}^{}=\Trace \R_{e^1}^{}\R_{e^2}^2=\frac{403}{8}, \ 
\Trace \R_{e^1}^2\R_{e^2}^2=\frac{197}{32}.
$$
Therefore, the dimensions of the eigenspaces are given by
$$
\EQ{ 
&d(0,0)=46851,\ 
d(0,{\textstyle \frac{1}{16}})=d({\textstyle\frac{1}{16}},0)=d({\textstyle \frac{1}{16},\frac{1}{16}})=47104,\cr
&d({\textstyle \frac{1}{2},0})=d({\textstyle 0,\frac{1}{2}})=2300,\ 
d({\textstyle \frac{1}{2},\textstyle \frac{1}{16}})=d({\textstyle \frac{1}{2},\textstyle \frac{1}{16}})=2048,\ 
d({\textstyle \frac{1}{2},\textstyle \frac{1}{2}})=23.\cr
}
$$
Therefore, for the idempotent $\idem=e^1+e^2$,
$$
\EQ{&d(0)=46851,\ 
d({\textstyle\frac{1}{16}})=94208,\ 
d({\textstyle\frac{1}{8}})=47104,\ 
d(\textstyle\frac{1}{2})=4600,\cr 
&
d(\textstyle\frac{9}{16})=4096,\ 
d(\textstyle1)=23,\ 
d(2)=2.\cr }
$$
In particular, the trace of the map
$$
\EQ{
\hfill 1
&\quad\hbox{on}\quad &B^\natural(0)
\oplus B^\natural({\textstyle\frac{1}{8}})
\oplus B^\natural({\textstyle\frac{1}{2}})
\oplus B^\natural({\textstyle1})
\oplus B^\natural({\textstyle2}),
\cr
-1
&\quad\hbox{on}\quad &B^\natural({\textstyle\frac{1}{16}})\oplus B^\natural({\textstyle\frac{9}{16}}),\cr
}
$$
which is the composition of two 2A involutions corresponding to
$e^1$ and $e^2$, is equal to
$46851-94208+47104+4600-4096+23+2=276$.  Hence this map is identified
with a 2B involution of the Monster.

We may do the same analysis for an embedding of $\Lvir(\frac12,0)^{\otimes
3}$.  However, the spectrum is not uniquely determined by the trace
formulae; there are two possibilities.
\subsubsection*{$\Lvir(\frac7{10},0)$ and 2A involution}
Let $U$ be the VOA $\Lvir(\frac7{10},0)$, for which the irreducible
modules are parametrized by the lowest conformal weight $h=0$, $3/80$,
$1/10$, $7/16$, $3/5$, $2/3$.  Hence the spectrum of the idempotent is
determined as
$$
\EQ{&d(0)=51054,\ 
d(\textstyle\frac{3}{80})=91392,\ 
d({\textstyle\frac{1}{10}})=47634,\ 
d(\textstyle\frac{7}{16})=4864,\ \cr
&d({\textstyle\frac{3}{5}})=1938,\ 
d(\textstyle\frac{3}{2})=1,
d(2)=1.\cr }
$$

The map
$$
\EQ{
\hfill 1
&\quad\hbox{on}\quad &B^\natural(0)\oplus B^\natural({\textstyle\frac{1}{10}})\oplus B^\natural({\textstyle\frac{3}{2}})\oplus B^\natural({\textstyle\frac{3}{5}})\oplus B^\natural(2),
\cr
-1
&\quad\hbox{on}\quad &B^\natural({\textstyle\frac{3}{80}})\oplus B^\natural({\textstyle\frac{7}{16}}),
\cr
}
$$
gives rise to an automorphism of $B^\natural$ by \cite{Miyamoto1}.
Since the trace is equal to $51054-91392+47634-4864+1938+1+1=4372$,
this map is identified with a 2A involution of the Monster.
\subsubsection*{$W_3$ algebra at $c=4/5$ and 3A element}
Let\footnote{The author presently does not know the existence of an
embedding of $W_3(\frac45)$ into $V^\natural$.} $U=W_3(\frac45)$ be the
vacuum sector of the $W_3$ algebra \cite{Zamolodchikov} at $c=4/5$.
It is isomorphic to $\Lvir(\frac45,0)\oplus \Lvir(\frac45,3)$ as
modules over the Virasoro algebra.  A realization of $W_3(\frac45)$ as
a VOA as well as its representation theory are described in \cite{KMY}
and \cite{Miyamoto4}.

There are 6 irreducible modules for this VOA, which are labeled as
$(h,\sigma)=(0,0)$, $(2/5,0)$, $({2}/{3},\pm)$, $({1}/{15},\pm)$,
where $h$ is the lowest conformal weight and $\sigma$ is the sign of
the eigenvalue of the action of a certain primary vector $w\in
W_3(\frac45)$ of conformal weight $3$.  Hence the spectrum is given by
$$
d(0)=57478,\ d(\textstyle\frac{1}{15})=129168,\
d({\textstyle\frac{2}{5}})=8671,\ d({\textstyle\frac{2}{3}})=1566,\
d(2)=1.
$$

Consider the case when the primary vector $w$ is also mapped to a
vector in the real form $V^\natural_\RR$. Then since
$[w_{(2)},t_{(1)}]=0$ and $w_{(2)}$ is alternating with respect to the
form $(\ |\ )$, we have a decomposition
$$
\EQ{
&B^\natural(\textstyle\frac{1}{15})
&=B^\natural(\textstyle\frac{1}{15},+)\oplus B^\natural(\textstyle\frac{1}{15},-),\quad B^\natural(\textstyle\frac{2}{3})
&=B^\natural(\textstyle\frac{2}{3},+)\oplus B^\natural(\textstyle\frac{2}{3},-),\cr
}
$$
into the sum of subspaces of equal dimensions for $h=1/15$ and $2/3$.
Now the map
$$
\EQ{
\hfill 1
&\quad\hbox{on}\quad &B^\natural(0)\oplus B^\natural({\textstyle\frac{2}{5}})\oplus B^\natural(2),\qquad
\zeta^{\pm1}
&\quad\hbox{on}\quad &B^\natural({\textstyle\frac{2}{3}},\pm)\oplus B^\natural({\textstyle\frac{1}{15}},\pm),\cr
}
$$
where $\zeta$ is a primitive 3rd root of unity, gives rise to an
automorphism of $B^\natural$ by \cite{Miyamoto4}. Since the trace is
equal to $57478+8671+(\zeta+\zeta^{-1})(129168+1566)/2+1=783$, it is
identified with a 3A element of the Monster.

This eigenspace decomposition is described in (24) of
\cite{MeyerNeutsch} and Lemma 4 of \cite{Norton}.  
\subsubsection*{$W_4$ algebra at $c=1$ and 4A element}
Let $U=W_4(1)$ be the vacuum sector of the $W_4$ algebra at $c=1$.  It
is realized as the fixed-point subspace $V^+_L$ of the lattice VOA
$V_L$ corresponding to the rank one lattice $L=\ZZ \gamma$ with
$\<\gamma,\gamma\>=6$ with respect to the $-1$ automorphism of the
lattice.  It is generated by the conformal vector $t$ and certain
primary vectors $w,z$ of conformal weight $3$ and $4$ respectively.
We may use the representation theory of $V^+_L$ developed in
\cite{DongNagatomo} and \cite{Abe}.

There are $10$ irreducible modules for this VOA, which are labeled as
$(h,\sigma)=(0,0)$, $(1,0)$, $({1}/{12},0)$, $(1/3,0)$, $(3/4,\pm)$, $(1/{16},\pm)$, $(9/{16},\pm)$.
Hence the spectrum is given by
$$
\EQ{
&d(0)=38226,\ 
d(\textstyle\frac{1}{16})=94208,\ 
d({\textstyle\frac{1}{12}})=48600,\
d({\textstyle\frac{1}{3}})=11178,\ \cr
&d(\textstyle\frac{9}{16})=4096,\ 
d(\textstyle\frac{3}{4})=552,\ 
d(1)=23,\ 
d(2)=1.\cr }
$$
These dimensions are determined as unique nonnegative integers
that satisfy the formula (\ref{Cor}), although the number of unknown
dimensions exceeds the number of equations.

Suppose that the primary vectors $w,z$ are mapped to the real form
$V^\natural_\RR$. Then the map
$$
\EQ{
\hfill 1
&\quad\hbox{on}\quad &B^\natural(0)\oplus B^\natural({\textstyle\frac{1}{3}})\oplus B^\natural(1)\oplus B^\natural(2),
\cr
\hfill-1
&\quad\hbox{on}\quad &B^\natural({\textstyle\frac{1}{12}})\oplus B^\natural({\textstyle\frac{3}{4}}),
\cr
\hfill\pm \sqrt{-1}
&\quad\hbox{on}\quad &B^\natural({\textstyle\frac{1}{16}},\pm)\oplus B^\natural({\textstyle\frac{9}{16}},\pm)
\cr
}
$$
gives rise to an automorphism of $B^\natural$ by the fusion rules of
$W_4(1)$. Since the trace is equal to
$38226-48600+11178-552+23+1=276$, it is identified with a 4A element
of the Monster.  Note that the trace of the square of this map, i.e.,
of the map
$$
\EQ{
\hfill 1
&\quad\hbox{on}\quad &B^\natural(0)\oplus
B^\natural({\textstyle\frac{1}{12}})\oplus
B^\natural({\textstyle\frac{1}{3}})\oplus
B^\natural(\textstyle\frac{3}{4})\oplus B^\natural(1)\oplus
B^\natural(2), \cr
-1 &\quad\hbox{on}\quad &B^\natural({\textstyle\frac{1}{16}})\oplus
B^\natural({\textstyle\frac{9}{16}}), \cr }
$$
is equal to $276$. Hence this map is identified with a 2B involution
of the Monster as expected.

This eigenspace decomposition is described in Lemma 5 of \cite{Norton}%
\footnote{There $24104+24104$ should read $47104+47104$.}.
\subsubsection*{$W_5$ algebra at $c=8/7$ and 5A element}
Unfortunately, the classification of irreducible modules and the
determination of fusion rules based on the theory of VOA for
$W_5(\frac{8}{7})$ seem to be missing.  However, formally applying the
expected properties of this algebra to our situation, the eigenspace
decomposition is expected to be
$$\EQ{
&\textstyle d(0)=27228,\ d(\frac2{35})=72010,\ d(\frac3{35})=76912,\ d(\frac{17}{35})=6688,\ d(\frac{23}{35})=1520,\cr
&\textstyle d(\frac27)=12122,\ d(\frac67)=133,\ d(\frac45)=268,\ d(\frac65)=2,\ d(2)=1.\cr}
$$
Since the trace of the map
$$
\EQ{
\hfill 1
&\quad\hbox{on}\quad &B^\natural(0)\oplus B^\natural({\textstyle\frac{2}{7}})\oplus B^\natural({\textstyle\frac{6}{7}})\oplus B^\natural(1)\oplus B^\natural(2),
\cr
\hfill \zeta^{\pm1}
&\quad\hbox{on}\quad &
B^\natural({\textstyle\frac{2}{35}},\pm)\oplus B^\natural({\textstyle\frac{17}{35}},\pm)\oplus B^\natural({\textstyle\frac{6}{5}},\pm),
\cr
\hfill \zeta^{\pm2}
&\quad\hbox{on}\quad &
B^\natural({\textstyle\frac{3}{35}},\pm)\oplus B^\natural({\textstyle\frac{23}{35}},\pm)\oplus B^\natural({\textstyle\frac{4}{5}},\pm),
\cr
}
$$
where $\zeta$ is a primitive 5th root of unity, is equal to 
$
27228+(\zeta+\zeta^{-1})(72010+6688+2)/2+(\zeta^2+\zeta^{-2})(76912+1520+268)/2+12122+133+1=134,
$
this map would be identified with a 5A element of the Monster (if we
appropriately choose the signs above\footnote{This ambiguity is a
matter of identification of the representations. There is no ambiguity
if we adopt the labeling as in \cite{FrenkelKacWakimoto}.}).

\Note \rm 
The fusion rules of $W_n$ algebras constructed by the quantized
Drinfeld-Sokolov reduction are determined by Frenkel et al.\
\cite{FrenkelKacWakimoto} via the Verlinde formula.  In particular,
the fusion ring of the first unitary series $W_n(c_n)$,
$c_n=2(n-1)/(n+2)$, is isomorphic to that of the level~$2$ integrable
highest weight representations of the affine Kac-Moody Lie algebra of
type $A_{n-1}^{(1)}$.  Then the map $[\lambda]\mapsto
\zeta_\lambda[\lambda]$, where
$\zeta_\lambda=\exp(2\pi\sqrt{-1}\sum_{i=1}^{n-1}im_i/n)$ for a
level~$2$ weight $\lambda=\sum_{i=1}^{n-1}m_i\barb{\Lambda}_i$, gives
an automorphism of the fusion ring over $\CC$.  Therefore, we expect
that an element of order $n$ of the Monster would be obtained by using
the decomposition of $V^\natural$ into the sum of irreducible
$W_n(c_n)$-modules.
\EndNote
\section{Generalization to Higher Degree}
Recall the notations and assumptions in section 2.3.  In particular,
$\Aut V$ is supposed to be finite. 
\subsection{The trace functions}
Let us set
\Equation
o(u)=u_{(n-1)}:V\rightarrow V
\EndEquation
for any $u\in V^n$ after Frenkel and Zhu \cite{FrenkelZhu}.  
Consider the trace functions
$$
\Trace o(a) q^{L_0}\quad\hbox{and}\quad \Trace o(a)o(b)q^{L_0}
$$
for elements $a,b\in B$, where $\Trace$ denotes the trace over the
whole space $V$.

In this subsection, we will express these trace functions, under the
corresponding assumptions, in terms of the character
\Equation
\character V=\sum_{n=0}^\infty (\dim V^n)q^n
\EndEquation
and the Eisenstein series $E_{2k}=E_{2k}(q)$, which we normalize as in
\cite{DongMason}:
\Equation
E_{2k}(q)=-\frac{B_{2k}}{(2k)!}+\frac{2}{(2k-1)!}\sum_{n=1}^\infty \left(\sum_{d|n}d^{2k-1}\right)q^n.
\EndEquation
Here $B_{m}$ are the Bernoulli numbers defined by
$$
\frac{t}{e^t-1}=\sum_{m=0}^\infty B_m\frac{t^m}{m!}.
$$
In particular,
\Equation
\EQ{
&E_2=-\frac{1}{12}+2q+6q^2+8q^3+14q^4+12q^5+\cdots\cr
&E_4=\frac{1}{720}+\frac{1}{3}q+3q^2+\frac{28}{3}q^3+\frac{73}{3}q^4+42q^5+\cdots\cr
}
\EndEquation

The result is summarized in the following theorem

\Theorem \sl 
Let $B$ be the Griess algebra of a VOA $V$ such that $\Aut V$ is finite.
\item
{\rm (1)}
If $V$ is of class~$\CS^2$ then, for any $a\in B$, 
$$
\Trace o(a)q^{L_0}
=\frac{2(a|\omega)}{c} q\frac{d}{dq} \character V.
$$
{\rm (2)}
If $V$ is of class~$\CS^4$ then, for any $a,b\in B$, 
$$
\EQ{
\Trace o(a)o(b)q^{L_0}
=&\left(\frac{44(a|b)+20(a|\omega)(b|\omega)}{c(5c+22)}\biggl(q\frac{d}{dq}\biggr)^2\right.\cr
&\ -(11+60E_2)\frac{c(a|b)-2 (a|\omega)(b|\omega)}{3c(5c+22)}q\frac{d}{dq}\cr
&\ \left.+(11+120E_2-720 E_4)\frac{c(a|b)-2(a|\omega)(b|\omega)}{360(5c+22)}\right)\character V.
\hfill\cr
}
$$
\label{Theorem3}
\EndTheorem

For instance,
$$
\EQ{
\Trace|_{V^3} &o(a)o(b)\cr
&=\frac{-2(20c^2+40 c\dim V^2+(3c-198)\dim V^3)}{c(5c+22)}(a|b)\cr
&+\frac{16(5c+10 \dim V^2+12 \dim V^3)}{c(5c+22)}(a|\omega)(b|\omega),\cr
\noalign{\vskip1ex}
\Trace|_{V^4} &o(a)o(b)\cr
&=\frac{-2(55c^2+98c \dim V^2+60c \dim V^3+(4c-352)\dim V^4)}{c(5c+22)}(a|b)\cr
&+\frac{4(55c+(5c+120) \dim V^2+60 \dim V^3+84\dim V^4)}{c(5c+22)}(a|\omega)(b|\omega).\cr
}
$$

We sketch the derivation of the formulae in the rest of this
subsection.  The formula~(1) immediately follows from (\ref{eq:a=})
and Lemma \ref{lemma:delta}: if $V$ is of class~$\CS^2$ then
$$
\Trace a_{(1)}q^{L_0}
=\frac{2(a|\omega)}{c}\Trace \omega_{(1)} q^{L_0}
=\frac{2(a|\omega)}{c}\Trace L_0 q^{L_0}
=\frac{2(a|\omega)}{c} q\frac{d}{dq} \character V.
$$
for any $a\in B$.

Suppose that $V$ is of class~$\CS^4$ and consider two vectors
$a,b\in B$ of the Griess algebra.  Then, by Proposition 4.3.5 of
\cite{Zhu}, we have 
\Equation 
\Trace o(a)o(b)q^{L_0}=\Trace
o(a_{[-1]}b)q^{L_0}-E_2\Trace o(a_{[1]}b)q^{L_0}-E_4\Trace
o(a_{[3]}b)q^{L_0}.
\label{eq:zhu}
\EndEquation
Here the operations $(a,b)\mapsto a_{[n]}b$, $(n\in\ZZ)$, are another
VOA structure on $V$ introduced by Zhu \cite{Zhu}, which is normalized
so that
\Equation
a_{[-1]}b
=a_{(-1)}b+\frac32 a_{(0)}b+\frac5{12}a_{(1)}b+\frac{11}{720}a_{(3)}b,\ 
a_{[1]}b
=a_{(1)}b-\frac16 a_{(3)}b,\ 
a_{[3]}b=a_{(3)}b.
\EndEquation
for $a,b\in B$.  Then, by (\ref{eq:zhu}) and (\ref{eq:projection}), we
have
\Equation
\EQ{
\Trace& o(a)o(b)q^{L_0}
=\frac{11}{720}\Trace o(\1)q^{L_0}+\frac{5(a|b)}{3c}\Trace o(\basis{2})q^{L_0}+\frac{3(a|b)}{c}\Trace o([3])q^{L_0}\cr
&+\frac{6c(a|b)-12(a|\omega)(b|\omega)}{c(5c+22)}\Trace o(\basis{4})q^{L_0}
+\frac{44(a|b)+20(a|\omega)(b|\omega)}{c(5c+22)}\Trace o(\basis{2,2})q^{L_0}\cr
&-E_2\frac{4(a|b)}{c}\Trace o(\basis{2})q^{L_0}
+\frac16 E_2(a|b)\Trace o(\1)q^{L_0}
-E_4(a|b)\Trace o(\1)q^{L_0}.\cr
}
\label{eq:zhucor}
\EndEquation
Hence we get Theorem \ref{Theorem3} (2) by using the following:
\vskip1ex
\vbox{
\Equation
\Trace o(\basis{2})=q\frac{d}{dq}\character V,\ 
\Trace o(\basis{3})=-2q\frac{d}{dq}\character V,\ 
\Trace o(\basis{4})=3q\frac{d}{dq}\character V,
\EndEquation
$$
\Trace o(\basis{2,2})q^{L_0}
=\left(\biggl(q\frac{d}{dq}\biggr)^2+\biggl(\frac{13}{6}+2E_2\biggr)q\frac{d}{dq}-\biggl(\frac{11c}{1440}+\frac{c}{12}E_2-\frac{c}{2}E_4\biggr)\right)\character V.
$$
}
\noindent
Here only the last one is not obvious.  By (\ref{eq:zhu}), we have
\Equation
\Trace o(\omega)o(\omega)q^{L_0}=\Trace o(\omega_{[-1]}\omega)q^{L_0}-E_2\Trace o(\omega_{[1]}\omega)q^{L_0}-E_4\Trace o(\omega_{[3]}\omega)q^{L_0}.
\label{eq:zhuomega}
\EndEquation
Substituting
$$
o(\omega_{[-1]}\omega)
=[2,2]_{(3)}-\frac{13}{6}[2]_{(1)}+\frac{11}{144}c,\ 
o(\omega_{[1]}\omega)=2[2]_{(1)}-\frac{c}{12},\ 
o(\omega_{[3]}\omega)=\frac{c}{2}.
$$
and
$$
\Trace o(\omega)q^{L_0}=\Trace L_0 q^{L_0}=q\frac{d}{dq} \character V,\ 
\Trace o(\omega)o(\omega)q^{L_0}=\biggl(q\frac{d}{dq}\biggr)^2\character V,
$$
we have the result.
\subsection{McKay-Thompson series for 2A involution}
Let $V$ be the moonshine module $V^\natural$.  In this subsection, we
will show that the McKay-Thompson series for a 2A involution is
determined by the formulae above using the fact that the character of
$V^\natural$ is given by
\Equation
\EQ{
\character V^\natural&=q(J(q)-744)\cr
&=1+196884 q^2+ 21493760q^3+864299970q^4+20245856256 q^5+\cdots\cr
}
\label{eq:Jfunction}
\EndEquation
and that the characters of $\Lvir(\frac12,h)$ are given by \Equation \EQ{
&\chi_0(q)&\textstyle=\character
\Lvir(\frac12,0)&=\frac12\left(\prod_{k=0}^\infty(1+q^{k+1/2})+\prod_{k=0}^\infty(1-q^{k+1/2})\right),\cr
&\chi_{1/2}(q)&\textstyle =\character
\Lvir(\frac12,\frac12)&=\frac12\left(\prod_{k=0}^\infty(1+q^{k+1/2})-\prod_{k=0}^\infty(1-q^{k+1/2})\right),\cr
&\chi_{1/16}(q)&\textstyle =\character \Lvir(\frac12,\frac1{16})&=q^{1/16}\prod_{k=1}^\infty
(1+q^k).\cr }
\label{eq:isingchar}
\EndEquation

Now, let $e$ be an idempotent of central charge $1/2$ in the
Conway-Griess algebra $B^\natural$ and consider the corresponding
Virasoro action $L^e_n=e_{(n+1)}$.  Consider the subspace
\Equation
P(h)=\{v\in V\,|\,L^e_nv=0\,\hbox{\ if\ $n\geq1$ \ and \
$L^e_0v=hv$}\}
\EndEquation
for each $h=0,1/2,1/16$, and set
\Equation
z_h(q)=\sum_{n=0}^\infty \dim (V^{n}\cap P(h))q^n.
\EndEquation
Then we have
\def\newdot#1{{\mathop{#1}\limits^{\textstyle\hskip3pt\bf.}}}
\def\newdots#1{{\mathop{#1}\limits^{\textstyle\hskip1.8pt\bf.\hskip-0.8pt.}}}
\Equation
\EQ{
\hfill z^{}_0(q)\chi^{}_0(q)+q^{-1/2}z^{}_{1/2}(q)\chi^{}_{1/2}(q)+q^{-1/16}z^{}_{1/16}(q)\chi^{}_{1/16}(q)&=\Trace q^{L_0},\cr
\hfill z^{}_0(q)\newdot{\chi}_0(q)+q^{-1/2}z^{}_{1/2}(q)\newdot{\chi}_{1/2}(q)+q^{-1/16}z^{}_{1/16}(q)\newdot{\chi}_{1/16}(q)&=\Trace o(e)q^{L_0},\cr
\hfill z^{}_0(q)\newdots{\chi}_0(q)+q^{-1/2}z^{}_{1/2}(q)\newdots{\chi}_{1/2}(q)+q^{-1/16}z^{}_{1/16}(q)\newdots{\chi}_{1/16}(q)&=\Trace o(e)^2q^{L_0},\cr
}
\label{eq:conditions}
\EndEquation
where $\newdot{\chi}(q)=q{d}/{dq}\chi(q)$ and
$\newdots{\chi}(q)=(q{d}/{dq})^2\chi(q)$.  By Theorem \ref{Theorem3},
\vskip1ex
\vbox{
\Equation
\Trace q^{L_0}=\character V^\natural,\quad\Trace o(e)q^{L_0}=\frac{1}{48}q\frac{d}{dq}\character V^\natural,
\EndEquation
$$\Trace o(e)^2q^{L_0}=\left(\frac{49}{13632}\biggl(q\frac{d}{dq}\biggr)^2
-\frac{47(11+60E_2)}{81792}q\frac{d}{dq}+\frac{47(11+120E_2-720 E_4)}{163584}\right)\character V^\natural,
$$
}
\noindent
where $\character V^\natural$ is given by (\ref{eq:Jfunction}).
Therefore, the condition (\ref{eq:conditions}) determines the series
$z_0(q)$, $z_{1/2}(q)$ and $z_{1/16}(q)$, so the McKay-Thompson series
$T_{\rm 2A}(q)$ via
\Equation
\EQ{
T_{\rm 2A}(q)
&\textstyle
=q^{-1}\Bigl(z^{}_0(q)\chi^{}_0(q)+q^{-1/2}z^{}_{1/2}(q)\chi^{}_{1/2}(q)-q^{-1/16}z^{}_{1/16}(q)\chi^{}_{1/16}(q)\Bigr).\cr
\cr}
\EndEquation
The result is written as a rational expression involving the functions
$J(q)$, $\chi_h(q)$, their first and the second derivatives and the
Eisenstein series $E_2(q)$ and $E_4(q)$.  We do not include the
explicit form in this paper.
\vfill\eject
\section*{Appendix}
\subsection*{A.1 \ The coefficients in the trace formula $(4)$}
{\small
$$\EQ{
&A_{1}=-c(2100 {c^5}+53650 {c^4}+304049 {c^3}-980942 {c^2}-1641936  
c+229152)\qquad\qquad\cr
&\hfill+(2455 {c^4}-193958 {c^3}+4032472 {c^2}+539488 c-1651584) d,\cr
&A_{2}=-c (1050 {c^5}+30965 {c^4}+279826 {c^3}+609848 {c^2}-271248  
c-150144)\cr
&\hfill-c (60 {c^4}-4929 {c^3}+96248 {c^2}+258428 c-56304) d,\cr
&A_{3}=
-c (1050 {c^5}+31085 {c^4}+270928 {c^3}+726848 {c^2}+1748472 c-79008)\cr
&\hfill+c(60 {c^4}-3969 {c^3}+20752 {c^2}+1761292 c+127440)d,\cr
&B=4c(1050{c^4}+30905{c^3}+289750{c^2}+281168c\cr
&\hfill-4d(120{c^4}-14853{c^3}+424928{c^2}+11132c-206448)),\cr
&C=8c(d-1)(120{c^3}-9437{c^2}+187858c+22968),\cr
&D=-192c(d-1)(100{c^2}-4297c-2852)(d-1),\quad E=15744c(30c+47).\cr
}
$$
}
\subsection*{A.2 \ The coefficients in the trace formula $(5)$}
{\small
$$\EQ{
&A_{1,2,3,4,5}= -5c(46200{c^6}+2154600{c^5}+31531073{c^4}\cr
&\hfill\hfill+123663366{c^3}-560461448{c^2}-1390398720c-168205824)\hfill\cr
&\hfill- 5d(100{c^6}-2405{c^5}-1037398{c^4}+70463896{c^3}-1249353984{c^2}+60544768c+766334976),\cr
& A_{1,2,4,3,5}\cr
&\quad=c(1500{c^5}-161985{c^4}+5500754{c^3}-19601928{c^2}-1338547904c-3497905152)(d-1),\cr
& A_{1,2,5,3,4}\cr
&\quad=c(-1500{c^5}+147745{c^4}-3380778{c^3}-83375368{c^2}+2968841472c+3711048192)(d-1),\cr
& A_{1,3,2,4,5}
= -c(115500{c^6}+5849050{c^5}+102135165{c^4}\cr
&\hfill\hfill+720684894{c^3}+1549368552{c^2}-664210624c+4461754368)\hfill\cr
&\hfill+ cd(300{c^5}-81505{c^4}+5253294\
{c^3}-87363968{c^2}-611758944c+4940713728),\cr
& A_{1,3,4,2,5}\cr
&\quad=c(500{c^5}-29035{c^4}+518574{c^3}-15730088{c^2}+553755136c-3442893312)(d-1),\cr
& A_{1,3,5,2,4}\cr
&\quad=c(-500{c^5}+14795{c^4}+1601402{c^3}-87247208{c^2}+1076538432c+3656036352)(d-1),\cr
& A_{1,4,2,3,5}
=-c(115500{c^6}+5848050{c^5}+102268115{c^4}\cr
&\hfill\hfill+715702714{c^3}+1553240392{c^2}+1228092416c+4516766208)\hfill\cr
&\hfill-cd(700{c^5}-51445{c^4}-271114{c^3}+83492128{c^2}-1280544096c-4995725568),\cr
& A_{1,4,3,2,5}\cr
&\quad=c(-500{c^5}+82955{c^4}-5023222{c^3}+122369272{c^2}-861258816c-1610972160)(d-1),\cr
& A_{1,4,5,2,3}\cr
&\quad=c(500{c^5}-118155{c^4}+6583582{c^3}-91119048{c^2}-815764608c+3601024512)(d-1),\cr
}$$
$$\EQ{
& A_{1,5,2,3,4}=-c(115500{c^6}+5847050{c^5}+102401065{c^4}\cr
&\hfill\hfill+710720534{c^3}+1557112232{c^2}+3120395456c+4571778048)\hfill\cr
&\hfill-cd(1700{c^5}-184395{c^4}+4711066{c^3}+79620288{c^2}-3172847136c-5050737408),\cr
& A_{1,5,3,2,4}\cr
&\quad=c(500{c^5}-49995{c^4}-41042{c^3}+118497432{c^2}-2753561856c-1665984000)(d-1),\cr
& A_{1,5,4,2,3}\cr
&\quad=c(-500{c^5}+103915{c^4}-4463606{c^3}-11858248{c^2}+2446058176c-3387881472)(d-1),\cr
& A_{2,3,1,4,5}\cr
&\quad=-3c(100{c^5}-21675{c^4}+907054{c^3}+11023128{c^2}-806389760c+1100745216)(d-1),\cr
& A_{2,4,1,3,5}\cr
&\quad=c(700{c^5}-67925{c^4}+2261018{c^3}-36941224{c^2}+526866240c-3357247488)(d-1),\cr
& A_{2,5,1,3,4}\cr
&\quad=c(1700{c^5}-200875{c^4}+7243198{c^3}-40813064{c^2}-1365436800c-3412259328)(d-1),\cr
&B_{1}=4c(115500{c^5}+5848250{c^4}+101927925{c^3}+740910478{c^2}+1067413032c+217343424)\cr
&\hfill + 4d(500{c^5}+288745{c^4}-25478878{c^3}+569319488{c^2}-269795104c-478959360),\cr
&B_{2}
=-4c(8100{c^4}-616655{c^3}+8745246{c^2}+142937384c-614801472)(d-1),\cr
&B_{3}
=4c(8100{c^4}-482575{c^3}-1572066{c^2}+339056296c-368532288)(d-1),\cr
&C=-8c(1780{c^4}-264997{c^3}+12872162{c^2}-203786696c-26642880)(d-1),\cr
&D=64c(3620{c^3}-510813{c^2}+15237868c+4458096)(d-1),\cr
&E=256c(2095{c^3}-161208{c^2}+3064358c+3847956)\
(d-1),\cr
&F=-3840c(3000{c^2}-125177c-223532)(d-1),\quad G=333312c(90c+259)(d-1),\cr
&H=-\frac{c}{12}(100{c^5}-13295{c^4}+498218{c^3}-387184\
{c^2}-189230304c-5501184)(d-1).\cr
}$$
}
\vskip2ex
\noindent
{\it Acknowledgement. }The author wishes to thank Alexander Ivanov,
Simon Norton and Michael Tuite for stimulating discussions and useful
comments on the subject.  He is deeply grateful to Ian Grojnowski for
hospitality at Cambridge.
\bibliographystyle{unsrt}

\end{document}